\documentclass[11pt]{article}%
\usepackage{amsmath}
\usepackage{amsfonts}
\usepackage{amssymb}
\usepackage{graphicx}
\providecommand{\U}[1]{\protect\rule{.1in}{.1in}}
\newcommand{\sfrac}[2]{\frac{\displaystyle #1}{\displaystyle #2}}
\newcommand{\ignore}[1]{}
\newcounter{subscript}

\begin{document}

\title{Algorithms for estimating spectral density functions for periodic potentials
on the half line}
\date{   }
\author{CHARLES FULTON \\
Department of Mathematical Sciences \\
Florida Insitute of Technology \\
Melbourne,  FL.  32901-6975 \\ \\
DAVID PEARSON \\
Department of Mathematics \\
University of Hull \\
Cottingham Road \\
Hull, HU6  7RX England \\ 
United Kingsom \\ \\
STEVEN PRUESS \\
1133 N. Desert Deer Pass \\
Green Valley,  Arizona  85614-5530
}
\maketitle


\abstract{For Hill's equation on $[0, \infty)$ we prove new characterizations
of the spectral function $\rho(\lambda)$ and the spectral density function
$f(\lambda)$ based on analysis involving a companion system of first order differential equations as in \cite{FPP2,FPP4}.
A numerical algorithm is derived and implemented based on
coefficient approximation.  Results for several examples, including the Mathieu
equation, are presented.}

\section{Introduction}

In this paper we consider Hill's equation, henceforth referred to as the
\textit{SL-equation},
\begin{equation}
\label{SL_EQ}-y^{\prime\prime}+ q(x) y = \lambda y, \qquad0 \le x < \infty,
\end{equation}
where $q(x)$ is real valued and periodic with period $\ell$, and we impose a boundary
condition
\begin{equation}
\label{SL_BC}y(0) \cos\alpha+ y^{\prime}(0) \sin\alpha= 0
\end{equation}
for some $\alpha\in[0, \pi)$.

Let a fundamental system of solutions for (\ref{SL_EQ}) be defined for all
$\lambda\in\mathrm{C\!\!\!I}$ by
\begin{equation}
\left[
\begin{array}
[c]{ll}%
\theta(0,\lambda) & \phi(0,\lambda)\\
\theta^{\prime}(0,\lambda) & \phi^{\prime}(0,\lambda)
\end{array}
\right]  =\left[
\begin{array}
[c]{cc}%
\cos\alpha & -\sin\alpha\\
\sin\alpha & \cos\alpha
\end{array}
\right]  ;\label{FSM}%
\end{equation}
we can define a unique Titchmarsh-Weyl $m-$function by (Im $\lambda\neq0$)
\begin{equation}
\theta(x,\lambda)+m(\lambda)\phi(x,\lambda)\in L_{2}(0,\infty
).\label{TITCH-WEYL}%
\end{equation}
The \textit{spectral function} is then defined for $\lambda \in [\Lambda,\infty)$ by the Titchmarsh-Kodaira
formula
\begin{equation}
\rho(\lambda)=\lim_{\epsilon\rightarrow0}\frac{\displaystyle1}%
{\displaystyle\pi}\int_{\Lambda}^{\lambda}\mbox{Im}\,m(\mu+i\epsilon
)\,d\mu,\label{TITCH-KOD}%
\end{equation}
where $\Lambda$ is the cutoff point for which equation (\ref{SL_EQ}) is
nonoscillatory in $(-\infty,\Lambda)$ and oscillatory in $(\Lambda,\infty)$,
or equivalently, the lowest point of the essential spectrum. The
\textit{spectral density function} is then defined for $\lambda \in [\Lambda,\infty)$ by
\begin{equation}
f(\lambda):=\rho^{\prime}(\lambda)=\frac{\displaystyle1}{\displaystyle\pi}%
\lim_{\epsilon\rightarrow0}\mbox{Im}[m(\lambda+i\epsilon)].\label{fDEF}%
\end{equation}

We now summarize some well known information on the spectrum associated with the
problem (\ref{SL_EQ})-(\ref{SL_BC}) (see, for example,  \cite[Chap 1-2]{EAST}).
For the case of periodic potentials with the above boundary condition, the
spectrum is known to be absolutely continuous and consisting of bands
interspersed with open intervals called gaps.  For the regular periodic problem
having the boundary conditions
\begin{equation}
y(0)=y(\ell),\qquad y^{\prime}(0)=y^{\prime}(\ell),\label{BC_PER}%
\end{equation}
let the eigenvalues be ordered by $\lambda_{0}\leq\lambda_{1}\leq\cdots$,
where eigenvalues of multiplicity two are written twice in the sequence. For
the regular semi-periodic problem having the boundary conditions
\begin{equation}
y(0)=-y(\ell),\qquad y^{\prime}(0)=-y^{\prime}(\ell),\label{BC_SEMI}%
\end{equation}
let the eigenvalues be ordered by $\mu_{0}\leq\mu_{1}\leq\cdots$, where
eigenvalues of multiplicity two are written twice in the sequence. Then the
eigenvalues of the periodic and semi-periodic problems occur in the order
\[
-\infty<\Lambda=\lambda_{0}<\mu_{0}\leq\mu_{1}<\lambda_{1}\leq\lambda_{2}%
<\mu_{2}\leq\mu_{3}<\lambda_{3}\leq\lambda_{4}\cdots  .
\]
If A denotes the self-adjoint operator associated with the SL problem  (\ref{SL_EQ})-(\ref{SL_BC})
then the closed intervals
\begin{equation}
\lbrack\lambda_{0},\mu_{0}],[\mu_{,}\lambda_{1}],[\lambda_{2},\mu_{2}%
],[\mu_{3},\lambda_{3}],\ldots\label{GAPS2}%
\end{equation}
constitute the essential spectrum $\sigma_{e}$ (or the stability set) of A. The complementary
set of open intervals
\begin{equation}
(\mu_{0},\mu_{1}),(\lambda_{1},\lambda_{2}),(\mu_{2},\mu_{3}),\ldots
\label{GAPS3}%
\end{equation}
are the gaps in $\sigma_{e}$, or the instability set. The spectral function
$\rho(\lambda)$ in (\ref{TITCH-KOD}) is absolutely continuous and monotone increasing on $\sigma_e$, so 
the absolutely continuous spectrum is $\sigma_{ac}(A)$ = $\sigma_e$.

In this paper we develop new characterizations for the spectral density
function $f(\lambda)$; one leads to a very efficient algorithm for its
calculation. In Section~2 we summarize some general information concerning the
SL equation (\ref{SL_EQ}) and a companion system of first order equations which we
will utilize in this paper. In Section~3 the new characterizations are derived. The remaining sections
develop the numerical scheme and show examples.

\section{Preliminaries}

\setcounter{equation}{0}

In this section we give the first order system of equations which we have
found to be a useful companion system for the study of the Sturm-Liouville
equation (the PQR-equations in \cite{FPP2},\cite{FPP4}),
introduce a standard basis for the solution space, and state some relations
which connect it to equation (\ref{SL_EQ}).

Consider the companion first order system for $U=(P,Q,R)^{T}$ for $\lambda
\in(\Lambda,\infty)$:%
\begin{equation}
\frac{\displaystyle dU}{\displaystyle dx}=\frac{d}{dx}\left[
\begin{array}
[c]{c}%
P\\
Q\\
R
\end{array}
\right]  =\left[
\begin{array}
[c]{ccc}%
0 & \lambda-q & 0\\
-2 & 0 & 2(\lambda-q)\\
0 & -1 & 0
\end{array}
\right]  \cdot\left[
\begin{array}
[c]{c}%
P\\
Q\\
R
\end{array}
\right].  \label{PQR_EQ}%
\end{equation}
.

The following statements are straightforward, if occasionally tedious, to verify.

\vskip6pt \noindent\textbf{1.} If $y$ is any solution of the SL-equation, then
$(\,(y^{\prime2},-2yy,y^{2})$ is a solution of equation (\ref{PQR_EQ}).

\vskip6pt \noindent\textbf{2}. If we let a fundamental system of the SL-equation be
defined by the initial conditions,%

\begin{equation}
\left[
\begin{array}
[c]{ll}%
u(0,\lambda) & v(0,\lambda)\\
u^{\prime}(0,\lambda) & v^{\prime}(0,\lambda)
\end{array}
\right]  =\left[
\begin{array}
[c]{cc}%
1 & 0\\
0 & 1
\end{array}
\right] ,\label{uv_EQ}%
\end{equation}
then a corresponding fundamental system of solutions of equation
(\ref{PQR_EQ}) is
\begin{equation}
U=\left[  U_{1},U_{2},U_{3}\right]  =\left[
\begin{array}
[c]{ccc}%
(u^{\prime})^{2} & u^{\prime}v^{\prime} & (v^{\prime})^{2}\\
-2uu^{\prime} & -[u^{\prime}v+uv^{\prime}] & -2vv^{\prime}\\
u^{2} & uv & v^{2}%
\end{array}
\right].  \label{U_BASIS}%
\end{equation}

\vskip6pt \noindent\textbf{3. \ }   If we represent a general solution of equation
(\ref{PQR_EQ}) in the form%
\begin{equation}
\left[
\begin{array}
[c]{c}%
P\\
Q\\
R
\end{array}
\right]  =a(\lambda)U_{1}(x,\lambda)+b(\lambda)U_{2}(x,\lambda)+c(\lambda
)U_{3}(x,\lambda), \label{U_COORD}%
\end{equation}
then (using the initial conditions (\ref{uv_EQ})) we have

\begin{equation}
a(\lambda)=R(0,\lambda),\qquad b(\lambda)=-Q(0,\lambda),\qquad c(\lambda
)=P(0,\lambda). \label{COORD}%
\end{equation}

\vskip6pt \noindent\textbf{4. \ } The solutions \ \{$\theta,\phi$\} \ defined
\ by the initial conditions (\ref{FSM}) are linearly related to the solutions
\{u, v\} (and vice versa) by
\begin{align}
\theta &  =u\cos\alpha+v\sin\alpha,\qquad\phi=-u\sin\alpha+v\cos
\alpha,\label{CB1}\\
u &  =\theta\cos\alpha-\phi\sin\alpha,\qquad v=\theta\sin\alpha+\phi\cos
\alpha.\label{CB2}%
\end{align}

\vskip6pt \noindent\textbf{5. \ }An indefinite inner product on the solution
space \ of equation \ (\ref{PQR_EQ}) \ may be defined by
\begin{equation}
\langle U_{1},U_{2}\rangle:=2(P_{1}R_{2}+P_{2}R_{1})-Q_{1}Q_{2} = const,\text{ independent of x }\in\lbrack0,\infty)%
\label{PQR_IP}%
\end{equation}
where $U_{k}=(P_{k},Q_{k},R_{k}),k=1,2$.

\vskip6pt \noindent\textbf{6.} \ For any solution $U$ = (P,Q,R)$^{T}$ of
equation (\ref{PQR_EQ}),%
\[
\frac{\displaystyle d}{\displaystyle dx}\langle U,U\rangle=\frac{\displaystyle
d}{\displaystyle dx}[4PR-Q^{2}]=0,
\]
i.e.%

\begin{equation}
4PR-Q^{2}=const,\text{ independent of x }\in\lbrack0,\infty)\label{PQR_NORM}%
\end{equation}

\vskip6pt \noindent\textbf{7.} \ If $U_{1}$ and $U_{2}$ are any two solutions
of equation (\ref{PQR_EQ}) represented as in (\ref{U_COORD}), then
\begin{equation}
\langle U_{1},U_{2}\rangle=2(a_{1}c_{2}+c_{1}a_{2})-b_{1}b_{2} \label{IP1}%
\end{equation}
and, in particular,
\begin{equation}
\langle U_{1},U_{1}\rangle=4a_{1}c_{1}-b_{1}^{2}. \label{IP2}%
\end{equation}

\vskip6pt \noindent\textbf{8.} \ If $U$ \ = (P,Q,R)$^{T}$ is any solution \ of
\ equation (\ref{PQR_EQ}) expressed as in (\ref{U_COORD}), then
\begin{equation}
\langle U,U\rangle=4PR-Q^{2}=4ac-b^{2}. \label{abc_INV1}%
\end{equation}

\vskip6pt \noindent\textbf{9.} If $U=(P,Q,R)^{T}$ is any solution of equation
(\ref{PQR_EQ}) expressed as in (\ref{U_COORD}), and if it is also written as%

\[
U=\gamma_{1}V_{1}+\gamma_{2}V_{2}+\gamma_{3}V_{3}%
\]
where
\begin{equation}
V=\left[  V_{1},V_{2},V_{3}\right]  =\left[
\begin{array}
[c]{ccc}%
(\theta^{\prime})^{2} & \theta^{\prime}\phi^{\prime} & (\phi^{\prime})^{2}\\
-2\theta\theta^{\prime} & -[\theta^{\prime}\phi+\theta\phi^{\prime}] &
-2\phi\phi^{\prime}\\
\theta^{2} & \theta\phi & \phi^{2}%
\end{array}
\right]  \label{abc_INV2}%
\end{equation}
is the fundamental system of (\ref{PQR_EQ}) generated by the solutions
\ \{$\theta,\phi$\} \ defined by \ (\ref{FSM}), then
\begin{equation}
4ac-b^{2}=4\gamma_{1}\gamma_{3}-\gamma_{2}^{2}.\label{abc_INV3}%
\end{equation}
In fact, the result holds if $\ \theta$ \ and $\ \phi$ \ are any two solutions
of the SL-equation with \ W$_{x}$($\theta(\cdot,\lambda),\phi(\cdot,\lambda)$)
= 1.

\vskip6pt \noindent\textbf{10. \ \ }If $y$ is any \ solution \ of \ the \ SL
\ equation  (\ref{SL_EQ})  and $U$=(P,Q,R)$^{T}$ is any solution of
\ companion system (\ref{PQR_EQ}) \ then
\[
\frac{\displaystyle d}{\displaystyle dx}[Py^{2}+Qyy^{\prime}+R(y^{\prime}%
)^{2}]=0,
\]
i.e.,%
\begin{equation}
P(x,\lambda)y^{2}(x,\lambda)+Q(x,\lambda)y(x,\lambda)y^{\prime}(x,\lambda
)+R(x,\lambda)(y^{\prime}(x,\lambda)^{2}=\text{ constant, independent of
\ x.}\label{PQR_INV}%
\end{equation}

\vskip6pt \ \ In the proofs in the next sections we will make frequent use of
the above results \ (particularly \textbf{1})\textbf{ } which relate the
solutions of the \ SL-equation \ to \ the solutions \ of the companion system
\ (\ref{PQR_EQ}). \ We\ exploited \ similar interrelations  in \ \cite{FPP2}
and \ \cite{FPP4} \ in the study of potentials on the half line satisfying
\ $q\in L_{1}(0,\infty)$.  Here \ we make use of the same interrelations in
the study of periodic potentials.\newline\noindent\textbf{Remark} : $\ $In our
previous papers  \cite{FPP2}, \cite{FPP4}  the system (\ref{PQR_EQ}) \ was
referred to as \ the \textquotedblleft PQR equations" \ (our notation);
\ however, \ the analysis leading to them (particularly \ the motivating
\ property \ (\ref{PQR_INV}) ) \ was discovered by M. Appell \ \cite{APPELL}
\ in \ 1880. Accordingly, we will henceforth refer to this first order system
\ as \ the \ Appell \ equations.

\section{\bigskip Characterizations of the spectral density function}

\setcounter{equation}{0}

In this section we give an analog of the closed form characterization
obtained in \cite{FPP4} when $q\in L_{1}(0,\infty)$. For the case of a
periodic potential on the half line $[0,\infty)$ the basic ideas from
\cite{AlN1}, \cite{AlN2}, \cite{FPP2}, and \cite{FPP4} carry over, at least for
values of $\lambda$ in the stability intervals, to yield several formulas for the
spectral density function.

We begin with the following definition as in \cite{AlN1}.

\vskip 6pt \noindent\textsc{Definition.} The Sturm-Liouville equation
(\ref{SL_EQ}), with $q(x)$ periodic of period $\ell$, satisfies
\textit{Condition A} for a given real value of $\lambda$ if and only if there
exists a complex-valued solution $y(x, \lambda)$ for which
\begin{equation}
\label{CondA}\lim_{N \to\infty} \frac{\displaystyle \int_{0}^{N} y(x,
\lambda)^{2} \, dx}{\displaystyle \int_{0}^{N} |y(x, \lambda)|^{2} \, dx} = 0.
\end{equation}

\vskip 6pt \noindent We now have the following lemmas. \vskip 6pt
\noindent\textsc{Lemma 1.} For $\lambda$ in the stability intervals, let
\begin{equation}
\label{psi1_DEF}\psi_{1}(x, \lambda) = p_{1}(x)\exp(ik(\lambda)x) 
\end{equation}
be the first Floquet solution for the characteristic exponent $\rho_{1} =
\exp(i \ell k(\lambda))$. Here the first Floquet solution for each $\lambda$ 
is understood to have the choice of $k(\lambda)$ such that  $0<\ell k(\lambda)<\pi$, and
$p_1(x)$  is periodic of period $\ell$. Then
\begin{equation}
\label{Eqn3-1}\lim_{N \to\infty} \frac{\displaystyle \int_{0}^{N} \psi_{1}(x,
\lambda)^{2} \, dx}{\displaystyle \int_{0}^{N} |\psi_{1}(x, \lambda)|^{2} \,
dx} = 0
\end{equation}
It follows from \cite[Theorem~2]{AlN1} that the spectrum of (\ref{SL_EQ}) with
(\ref{SL_BC}) is absolutely continuous in the stability intervals.

\vskip 6pt \noindent\textsc{Proof.}
\begin{align*}
\int_{0}^{m \ell} \psi_{1}^{2} \, dx  &  = [1 + \sum_{j=1}^{m-1} \exp(2i j \ell k(\lambda) )]
\int_{0}^{\ell} \exp(2ik(\lambda)t) p_{1}^{2}(t) \, dt\\
&  = \frac{\displaystyle 1 - [\exp(2i \ell k(\lambda))]^{m}}{\displaystyle 1 - \exp(2i\ell k(\lambda))} \int_{0}^{\ell} \exp(2ik(\lambda)t) p_{1}^{2}(t) \, dt.
\end{align*}
But $0 < k < \pi/ \ell$, so $\exp(2i \ell k(\lambda)) \ne1$ for $\lambda$ in the
stability intervals; hence, the ratio is bounded above. It follows that
\[
\frac{\displaystyle \left|  \int_{0}^{m \ell} \psi_{1}(x, \lambda)^{2} \, dx
\right|  }{\displaystyle \int_{0}^{m \ell} |\psi_{1}(x, \lambda)|^{2} \, dx}
\le\frac{\displaystyle K \int_{0}^{\ell} |p_{1}|^{2} \, dt}%
{\displaystyle \int_{0}^{m \ell} |p_{1}|^{2} \, dt} \to0
\]
as $m \to\infty$. \quad\rule{2.2mm}{3.2mm}

\vskip 6pt \noindent\textsc{Lemma 2.} Let $\lambda\in\cup_{m=0}^{\infty}
[\lambda_{2m}, \mu_{2m}] \cup[\mu_{2m+1}, \lambda_{2m+1}]$. Then, since
Condition~A holds for $\lambda$ in these stability intervals, there is a
complex-valued function $\xi(\lambda)$ that is uniquely defined for $\lambda$
in the stability intervals by the properties

(i) $\mbox{Im } \xi(\lambda) > 0$, \noindent and

(ii)
\[
\lim_{N \to\infty} \frac{\displaystyle \int_{0}^{N} \left( \theta(x, \lambda) +
\xi(\lambda) \phi(x, \lambda) \right)^{2} \, dx}{\displaystyle \int_{0}^{N}
|\theta(x, \lambda) + \xi(\lambda) \phi(x, \lambda)|^{2} \, dx} = 0.
\]

\vskip6pt \noindent\textsc{Proof.} This is proved in \cite[Lemma1]{AlN1}.
\quad\rule{2.2mm}{3.2mm}

\vskip 6pt Since the solution that satisfies Condition~A is unique up to a
constant multiple, it follows from Lemma~1 and Lemma~2 that there exists a
constant $K \ne0$ such that
\begin{equation}
\label{Eqn3-3}\psi_{1}(x, \lambda) = K[ \theta(x, \lambda) + \xi(\lambda)
\phi(x, \lambda) ].
\end{equation}
From (\ref{psi1_DEF}) at $x = 0$ we have
\[
\psi_{1}(0) \cos\alpha+ \psi_{1}^{\prime}(0) \sin\alpha= p_{1}(0) \cos\alpha+
[ik p_{1}(0) + p_{1}^{\prime}(0)] \sin\alpha.
\]
But from (\ref{Eqn3-3})
\begin{align*}
\psi_{1} (0) \cos\alpha+ \psi_{1}^{\prime}(0) \sin\alpha &  = K[ \theta(0,
\lambda) + \xi(\lambda) \phi(0, \lambda) ] \cos\alpha+ K[ \theta^{\prime}(0,
\lambda) + \xi(\lambda) \phi^{\prime}(0, \lambda) ] \sin\alpha\\
&  = K[\cos^{2} \alpha- \xi\sin\alpha\cos\alpha+ \sin^{2} \alpha+ \xi
\sin\alpha\cos\alpha]\\
&  = K.
\end{align*}

It also follows immediately from Theorem~2 of \cite{AlN1} that for all
$\lambda$ in the stability intervals, we have that the function $\xi(\lambda)$
is, in fact, the boundary value of the Titchmarsh-Weyl $m-$function defined in
(\ref{TITCH-WEYL}), that is,
\begin{equation}
\xi(\lambda)=A(\lambda)+iB(\lambda):=\lim_{\epsilon\rightarrow0}%
m(\lambda+i\epsilon). \label{Eqn3-5}%
\end{equation}
We are now ready to prove the following theorem:

\vskip6pt \noindent\textsc{Theorem~1.} For $\lambda$ in the stability
intervals, there exists a solution $U=(P,Q,R)^{T}$ of Appell's equation
\ (\ref{PQR_EQ}), unique up to a constant multiple, which is periodic of
period $\ell$ on $(0,\infty)$.

\vskip6pt \noindent\textsc{Proof.} For $\lambda$ in the stability intervals we
set
\begin{equation}
U:=\left[
\begin{array}
[c]{c}%
P\\
Q\\
R
\end{array}
\right]:=\left[
\begin{array}
[c]{c}%
\psi_{1}^{\prime}\psi_{2}^{\prime}\\
-(\psi_{1}\psi_{2}^{\prime}+\psi_{2}\psi_{1}^{\prime})\\
\psi_{1}\psi_{2}%
\end{array}
\right]  ,\label{T1-1}%
\end{equation}
where $\psi_{1}(x,\lambda)=p_{1}(x)\exp(ik(\lambda)x)$, and $\psi
_{2}={\overline{\psi}_{1}}$. Since $p_{1}(x)$ is periodic of period $\ell$, it
follows that $p_{1}^{\prime}$ is also periodic of period $\ell$. We obtain
from (\ref{T1-1}) that
\begin{equation}
\left[
\begin{array}
[c]{c}%
P\\
Q\\
R
\end{array}
\right]  =\left[
\begin{array}
[c]{c}%
k(\lambda)^{2}|p_{1}|^{2}+|p_{1}^{\prime}|^{2}-2k(\lambda)\mbox{ Im}(p_{1}\overline{p_{1}^{\prime}})\\
-2\mbox{Re}(p_{1} \overline{p_{1}^{\prime}})\\
|p_{1}|^{2}%
\end{array}
\right]  \label{T1-2}%
\end{equation}
and each component is real valued and periodic with period $\ell$.

To prove that the periodic solution is unique, up to constant multiple,
consider the fundamental solution matrix of Appell's \ system of equations obtained by
replacing \ \{u,v\} in \ (\ref{U_BASIS}) \ by the Floquet solutions
\{$\psi_{1},\psi_{2}$\}, and let T(x) \ be the transfer \ matrix
\ which carries \ $U(x,\lambda)$ \ to \ $U(x+\ell,\lambda)$, \ i.e.%
\[
U(x+\ell,\lambda)=T(x)U(x,\lambda).
\]
Since the Floquet solutions satisfy \ $\psi_{j}(x+\ell,\lambda)=\rho_{j}%
\psi_{j}(x,\lambda)$, \ where \ $\rho_{j}=\exp(\pm ik(\lambda)\ell),\text{ } j=1,2,$
\ are the Floquet exponents for \ $\lambda$ \ in \ the stability intervals, it
follows that the first and third columns of (\ref{U_BASIS}) \ (with $\psi_{1}$
and $\psi_{2}$) are eigenvectors of T(x) with eigenvalues \ $\rho_{1}^{2}$ and
\ $\rho_{2}^{2}$ \ (which are not one), and the second \ column \ is an
eigenvector of \ T(x) \ with \ eigenvalue \ $\rho_{1}\rho_{2}$ $=1.\ $Hence T
\ has  a one-dimensional eigenspace \ for \ which P, Q and R are  all
periodic of period \ $\ell$. \quad\rule{2.2mm}{3.2mm}


\vskip6pt \noindent The next result provides three different representations
for the spectral density function $f(\lambda)$.

\vskip6pt \noindent\textsc{Theorem 2}. For $\lambda$ in a stability interval
\ let \ $U=(P,Q,R)^{T}$ \ be the periodic solution of Appell's system
\ (\ref{PQR_EQ}) which is normalized by (compare (\ref{PQR_NORM}))%

\begin{equation}
\langle U,U\rangle=4PR-Q^{2}=4.\label{U_NORM}
\end{equation}
Let \ \{$a(\lambda),b\left(  \lambda\right)  ,c\left(  \lambda\right)  $\}
\ be the coefficients in the representation \ (\ref{U_COORD}) \ of this
periodic solution. Then the spectral density function defined by (\ref{fDEF})
admits the following representations:\bigskip%
\begin{align}
f(\lambda) &  =\left\vert \frac{\displaystyle1}{\displaystyle\pi\lbrack
c(\lambda)\sin^{2}\alpha+b(\lambda)\sin\alpha\cos\alpha+a(\lambda)\cos
^{2}\alpha]}\right\vert \label{f1}\\
&  =\left\vert \frac{\displaystyle1}{\displaystyle\pi\lbrack P(0,\lambda
)\sin^{2}\alpha-Q(0,\lambda)\sin\alpha\cos\alpha+R(0,\lambda)\cos^{2}\alpha
]}\right\vert \label{f2}\\
&  =\left\vert \frac{\displaystyle1}{\displaystyle\pi\lbrack P(x,\lambda
)\phi(x,\lambda)^{2}+Q(x,\lambda)\phi(x,\lambda)\phi^{\prime}(x,\lambda
)+R(x,\lambda)\phi^{\prime}(x,\lambda)^{2}]}\right\vert .\label{f3}%
\end{align}
Here it will be observed that the normalization (\ref{U_NORM}) fixes the
periodic solution only up to a $\ \pm$ \ sign; it is for this reason that we
take the absolute value sign in these formulas to ensure that \ $f(\lambda)$
$\geq0,$ as required.  In the applications it often happens that the denominators in the
above expressions are positive in one stability interval and negative in another.

\vskip6pt \noindent\textsc{Proof}. \ The formulas \ (\ref{f1}) and (\ref{f2})
\ are equivalent because the representation (\ref{U_COORD}) \ guarantees that
\ \{$a(\lambda),b(\lambda),c(\lambda)$\} are given by (\ref{COORD}). \ The
denominator in (\ref{f3}) \ is constant, independent of x, by (\ref{PQR_INV})
\ and equal to (\ref{f2}) on evaluation at \ $x=0$. \ So it suffices to prove
(\ref{f2}) subject to the normalization \ (\ref{U_NORM}). \ Since $\psi
_{1}(x,\lambda)$ \ is linearly dependent on \ $\theta(x,\lambda)+\xi
(\lambda)\phi(x,\lambda)$ \ by (\ref{Eqn3-3}) \ (where $\ \xi(\lambda)$ is the
complex valued function defined on the stability intervals in Lemma 2), and
$\psi_{2}=\overline{\psi_{1}}$ \ is linearly dependent on \ $\theta(x,\lambda)+\overline{\xi\left(
\lambda\right)}  \phi(x,\lambda)$, \ we may represent the periodic solution in
(\ref{T1-1}) as
\begin{equation}
U(x,\lambda):=\left[
\begin{array}
[c]{c}%
P(x)\\
Q(x)\\
R(x)
\end{array}
\right]:=K(\lambda)\left[
\begin{array}
[c]{c}%
(\theta^{\prime}+\xi\phi^{\prime})(\theta^{\prime}+{\overline{\xi}}%
\phi^{\prime})\\
-[(\theta+\xi\phi)(\theta^{\prime}+{\overline{\xi}}\phi^{\prime}%
)+(\theta+{\overline{\xi}}\phi)(\theta^{\prime}+\xi\phi^{\prime})]\\
(\theta+\xi\phi)(\theta+{\overline{\xi}}\phi)
\end{array}
\right] \label{T2-1}
\end{equation}
for some real constant \ $K(\lambda)$, independent of x. \ The required
nomalization \ (\ref{U_NORM}) \ is equivalent by \ (\ref{abc_INV1}) to%

\begin{equation}
4a(\lambda)c(\lambda)-(b(\lambda))^{2}=4. \label{T2-2}
\end{equation}
Using the representation of the periodic solution in terms of the fundamental
system (\ref{abc_INV2}) \ of Appell's equations,%
\begin{equation}
\left(
\begin{array}
[c]{c}%
P\\
Q\\
R
\end{array}
\right)  =\gamma_{1}V_{1}+\gamma_{2}V_{2}+\gamma_{3}V_{3},\label{V_COORD}
\end{equation}
and comparing the \ $R$-component with the $R$-component in (\ref{T2-1}), gives\ %

\[
\gamma_{1}=K(\lambda),\;\gamma_{2}=2\operatorname{Re}(\xi(\lambda
))K(\lambda),\text{ \ and}\;\gamma_{3}=\left\vert \xi(\lambda)\right\vert
^{2}K(\lambda).
\]
Hence from \ (\ref{abc_INV1}) \ and \ (\ref{abc_INV3}) \ we have%

\begin{align*}
4PR-Q^{2} &  =4ac-b^{2}\\
&  =4\gamma_{1}\gamma_{3}-\gamma_{2}^{2}\\
&  =4(K(\lambda))^{2}\left[  \left\vert \xi(\lambda)\right\vert ^{2}%
-\operatorname{Re}^{2}\xi(\lambda)\right]  \\
&  =4(K(\lambda))^{2}\left[  \operatorname{Im}^{2}\xi(\lambda)\right]  =4
\end{align*}
if and only if%

\begin{equation}
\lbrack K(\lambda)\operatorname{Im}\xi(\lambda)]^{2}=1.\label{T2-3}
\end{equation}

From (\ref{Eqn3-5}) and \ (\ref{fDEF}) \ it follows (since
\ $\operatorname{Im}\xi(\lambda)>0$ by Lemma 2) that \ the normalization
(\ref{U_NORM}) holds if and only if%

\begin{equation}
f(\lambda):=\frac{1}{\pi}\operatorname{Im}\xi(\lambda)=\frac{1}{\pi\left\vert
K(\lambda)\right\vert }.\label{T2-4}
\end{equation}
Next, we use the initial conditions (\ref{FSM}) \ to evaluate the right hand
side of (\ref{T2-1}) and then substitute into the denominator of (\ref{f2}) to obtain,%

\begin{align}
&  \pi\left[  P(0,\lambda)\sin^{2}\alpha-Q(0,\lambda)\sin\alpha\cos
\alpha+R(0,\lambda)\cos^{2}\alpha\right]\nonumber  \\
&  =\pi K(\lambda)\left[
\begin{array}
[c]{c}%
(\sin^{2}\alpha+\cos^{2}\alpha)\cdot1\\
+(\sin^{3}\alpha\cos\alpha+\sin\alpha\cos^{3}\alpha-\sin^{3}\alpha\cos
\alpha-\sin\alpha\text{ }\cos^{3}\alpha)(\xi+\overline{\xi})\\
+(2\sin^{2}\alpha\text{ }\cos^{2}\alpha-2\sin^{2}\alpha\text{ }\cos^{2}%
\alpha)\left\vert \xi\right\vert ^{2}%
\end{array}
\right] \nonumber \\
&  =\pi K(\lambda)\label{T2-5}
\end{align}
The formula (\ref{f2}) \ now follows \ from \ (\ref{T2-4}) \ and
\ (\ref{T2-5}). \quad\rule{2.2mm}{3.2mm}

\bigskip

\vskip 6pt \noindent A more useful characterization of $f(\lambda)$ is given
by the following result.

\vskip6pt \noindent\textsc{Theorem 3}. Assume $\lambda$ is in a stability
interval. Then
\begin{equation}
f(\lambda)=\left\vert \frac{\displaystyle\sqrt{4-[u(\ell,\lambda)+v^{\prime
}(\ell,\lambda)]^{2}}}{\displaystyle2\pi\left[  u^{\prime}(\ell,\lambda
)\sin^{2}\alpha+(u(\ell,\lambda)-v^{\prime}(\ell,\lambda))\cos\alpha\sin
\alpha-v(\ell,\lambda)\cos^{2}\alpha\right]  }\right\vert .\label{f4}%
\end{equation}
Here the absolute value is needed to ensure that $f(\lambda)\geq0$; 
this is due to the fact that the denominator could be
negative in some of the stability intervals, and also corresponds to the fact that
the normalization of \{$a,b,c$\} in (\ref{NORM}) fixes \{$a,b,c$\} only up to a $\pm$ sign.

\vskip6pt \noindent\textsc{Proof.} From Theorem 2 it follows that if we can
construct a solution \ $U=(P,Q,R)^{T}$ of Appell's first order system
(\ref{PQR_EQ}) \ which is periodic of period $\ell$ and satisfies the
normalization \ $4PR-Q^{2}=4$ then we can use it to get $f(\lambda)$ (e.g.
from \ any one of the formulas \ (\ref{f1}), (\ref{f2}), or (\ref{f3})). In
particular, if this periodic solution is represented in the form
(\ref{U_COORD}) it follows from (\ref{abc_INV1}) that the coefficients
\{$a,b,c$\} satisfy%

\begin{equation}
4a(\lambda)c(\lambda)-b^{2}(\lambda)=4P(x,\lambda)R(x,\lambda)-Q(x,\lambda
)^{2}=4. \label{NORM}
\end{equation}
Considering only the third component of \ (\ref{U_COORD}) it therefore
suffices to generate coefficients \ \{$a,b,c$\} for which the quadratic form%

\begin{equation}
R(x)=a(u(x,\lambda))^{2}+bu(x,\lambda)v(x,\lambda)+c(v(x,\lambda))^{2} \label{QUAD_FORM}%
\end{equation}
is periodic of period $\ell$ and such \ that (\ref{NORM}) holds for $\lambda$ \ in
the stability intervals. \ Then \ $f(\lambda)$ is given by (\ref{f1}) with
this choice of \ \{$a,b,c$\}; or by (\ref{f2}), (\ref{f3}) \ where
\ \{$P,Q,R$\} is the corresponding periodic solution \ (\ref{U_COORD}) \ of
\ Appell's equations. 
To manufacture \ \{$a,b,c$\} consider \ the \ SL-equation (\ref{SL_EQ}) in the
system form%

\begin{align*}
\frac{d}{dx}\Psi & =\left(
\begin{array}
[c]{ll}%
0 & 1\\
-(\lambda-q(x)) & 0
\end{array}
\right)  \Psi =A(x)\cdot\Psi,\qquad \Psi(x,\lambda)=\left(
\begin{array}
[c]{l}%
y\\
y^{\prime}%
\end{array}
\right) .
\end{align*} 
Let $\ \Psi_{\alpha}(\cdot,\lambda)$ \ be the solution of the initial value problem%

\begin{equation}
\Psi_{\alpha}^{\prime}(x)=A(x)\Psi_{\alpha}(x),\qquad\Psi_{\alpha}(\alpha
)=I,\;\alpha\in\lbrack0,\infty).\label{IVP}
\end{equation}
The following facts are easily verified:

\begin{align}
\Psi_{0}(x)  & =\left(
\begin{array}
[c]{ll}%
u(x) & v(x)\\
u^{\prime}(x) & v^{\prime}(x)
\end{array}
\right) \label{T3-1} \\
\left[  \Psi_{0}(x)\right]  ^{-1}  & =\Psi_{x}(0)\label{T3-2}\\
\forall \text{ } \text{solutions}\;y\;of\;(\ref{SL_EQ})  & :\left(
\begin{array}
[c]{l}%
y(x)\\
y^{\prime}(x)
\end{array}
\right)  =\Psi_{0}(x)\left(
\begin{array}
[c]{l}%
y(0)\\
y^{\prime}(0)
\end{array}
\right) \label{T3-3} \\
\forall \text{ } \text{solutions}\;y\;of\;(\ref{SL_EQ})  & :\left(
\begin{array}
[c]{l}%
y(0)\\
y^{\prime}(0)
\end{array}
\right)  =\Psi_{x}(0)\left(
\begin{array}
[c]{l}%
y(x)\\
y^{\prime}(x)
\end{array}
\right)\label{T3-4}  \\
\Psi_{\alpha}(x)  & =\Psi_{t}(x)\cdot\Psi_{\alpha}(t) \label{T3-5}
\end{align}
The fact that \ $q(x)$ has period \ $\ell$ \ implies%

\begin{equation}
\Psi_{\ell}(x+\ell)=\Psi_{0}(x),\label{T3-6}
\end{equation}
and hence that $\Psi_{x}(x+\ell)$ \ is periodic with period $\ell.$ To
generate a quadratic form in \ $u$ and \ $v$ \ which is periodic of period
\ $\ell$ we put%

\begin{align}
\Phi(x)  & :=\left(
\begin{array}
[c]{ll}%
1 & 0
\end{array}
\right)  \Psi_{x}(x+\ell)\left(
\begin{array}
[c]{l}%
0\\
1
\end{array}
\right) \nonumber \\
& =\left(
\begin{array}
[c]{ll}%
1 & 0
\end{array}
\right)  \Psi_{\ell}(x+\ell)\Psi_{0}(\ell)\Psi_{x}(0)\left(
\begin{array}
[c]{l}%
0\\
1
\end{array}
\right) \nonumber \\
& =\left(
\begin{array}
[c]{ll}%
1 & 0
\end{array}
\right)  \Psi_{0}(x)\Psi_{0}(\ell)\Psi_{0}(x)^{-1}\left(
\begin{array}
[c]{l}%
0\\
1
\end{array}
\right)\nonumber  \\  \\
& =\left(
\begin{array}
[c]{ll}%
u(x) & v(x)
\end{array}
\right)  \left[
\begin{array}
[c]{ll}%
u(\ell) & v(\ell)\\
u^{\prime}(\ell) & v^{\prime}(\ell)
\end{array}
\right]  \left(
\begin{array}
[c]{l}%
-v(x)\\
u(x)
\end{array}
\right) \nonumber \\ \\
& =v(\ell)(u(x))^{2}-\left[  u(\ell)-v^{\prime}(\ell)\right]
u(x)v(x)-u^{\prime}(\ell)(v(x))^{2}. \label{T3-7}
\end{align}
Since $\ \Phi(x+\ell)=\Phi(x)$ \ for all \ $x\in\lbrack0,\infty)$, we need
\ only normalize the coefficients to achieve the required normalization (\ref{NORM}).
\ Taking $\ \gamma\cdot\Phi(x)$ so that \ $4ac-b^{2}=4,$ we find%

\begin{align}
& \gamma^{2}\left[  -4v(\ell)u^{\prime}(\ell)-(u(\ell)-v^{\prime}(\ell))^{2}\right]
\\
& =\gamma^{2}\left[  4-(u(\ell)+v^{\prime}(\ell))^{2}\right]  \\
& =4,
\end{align}
so that the required \ normalization is achieved with%

\begin{equation}%
\left(
\begin{array}
[c]{l}%
a(\lambda)\\
b(\lambda)\\
c(\lambda)
\end{array}
\right)
=\frac{1}{\sqrt{4-(u(\ell)+v^{\prime}(\ell))^{2}}}\left(
\begin{array}
[c]{l}%
-2v(\ell)\\
2(u(\ell)-v^{\prime}(\ell))\\
2u^{\prime}(\ell)
\end{array}
\right) \label{T3-8} ,
\end{equation}
and substitution of this into (\ref{f1}) yields \ (\ref{f4}). \quad\rule{2.2mm}{3.2mm}

\vskip 18pt

\section{The Numerical Method}                          
\setcounter{equation}{0}

In this section and the following two sections we describe a new numerical algorithm for obtaining approximations to the spectral density function, by making use of the representation (\ref{f4}) in Theorem 3, and compare performance with {\bf SLEDGE}. For general information and discussion of numerical methods for Sturm-Liouville problems we refer to Pryce's book \cite{PRYCE}, and for the the computation of spectral functions using the method of {\bf SLEDGE} we refer to our previous papers \cite{SLEDGE,TOMS,TOMSEAST}.  In contrast to {\bf SLEDGE}, the above  Theorem 3 for periodic potentials  enables computation of the spectral density function on the stability intervals by shooting (with piecewise trigonometric / hyperbolic splines) over a single period.

To compute $u$ and $v$ we employ the method of coefficient
approximation by which $q$ is replaced by a step-function
approximation $\hat q$.  We write the analog to (\ref{SL_EQ}) as
\begin{equation}    \label{yHat_EQ}
 - \hat y'' + \hat q(x) \hat y = \lambda \hat y, \qquad a \le x < \infty
\end{equation}
and $\hat u$ and $\hat v$ will satisfy (\ref{yHat_EQ}) with initial
conditions analogous to those of $u$ and $v$, respectively.  We
extend the formula (\ref{f4}) by defining for any $\lambda$
\begin{equation}  \label{fHat}
   \hat f(\lambda) = \sfrac{\sqrt{\max\{0, 4 - [\hat u(\ell) + 
        \hat v'(\ell)]^2\}}} {2 \pi |\hat u'(\ell) \sin^2 \alpha
        + (\hat u(\ell) - \hat v'(\ell)) \sin \alpha \cos \alpha - 
           \hat v(\ell) \cos^2 \alpha|}
\end{equation}
as an estimate of $f$.  The appeal of this approach is that closed-form 
solutions, piecewise circular or hyperbolic trig functions, are known 
for $\hat u$ and $\hat v$, admitting efficiencies of computation and 
analysis.

There are two potential numerical challenges in trying to integrate
(\ref{yHat_EQ}) and use (\ref{fHat}):
(1) mathematical instability when $\lambda < q(x)$, and 
(2) loss of accuracy if the numerator and demominator of (\ref{fHat}) 
vanish simultaneously.  We note that these difficulties arise in the 
original equations (\ref{SL_EQ}) and (\ref{f4}), so we would expect them 
to be inherited by any computational approach.  With coefficient 
approximation it is straightforward to address both of these issues.

In \cite{SLEDGE} we presented a stabilizing algorithm to solve 
(\ref{yHat_EQ}) for the regular Sturm-Liouville problem, which is a 
boundary value problem.  A similar approach will work here.  First,
we provide more detail of the algorithm.  We first subdivide $[0, \ell]$
into $N$ intervals
\[  0 = x_1 < x_2 < \cdots  < x_{N+1} = \ell; \]
we set $h_n = x_{n+1} - x_n$to be the width of the $n$th subinterval. 

On any subinterval $(x_n, x_{n+1})$ we choose $\hat q(x) = q_n$ to be 
constant (usually the $q$ value at the midpoint); then the differential 
equation (\ref{yHat_EQ}) has the closed-form solution 
\begin{equation}      \label{yHat}
 \hat y(x) = \hat y(x_n) \phi_n' (x-x_n) + \hat y'(x_n)
             \phi_n (x-x_n)
\end{equation}
with
\[   \phi_n (t)  =  \left\{ \begin{array}{ll}
                    \sin \omega_n t / \omega_n &  \tau_n  > 0 \\
                    \sinh \omega_n t / \omega_n  &  \tau_n  < 0 \\
                     t   & \tau_n = 0,
                    \end{array} \right.      \]
where
\begin{equation}    \label{TAU}
  \tau_n = \lambda - q_n 
\end{equation}
and
\[  \omega_n = \sqrt{|\tau_n|} .   \]
It follows that
\begin{equation}      \label{dyHat}
 \hat y'(x) = - \tau_n \hat y(x_n) \phi_n (x-x_n) + \hat y'(x_n)
             \phi_n' (x-x_n).
\end{equation}

In practice, one should use a truncated series expansion for small 
$|\tau_n| h_n^2$, e.g., 
\[   \phi_n (t) = t[1-\tau_n t^2/6 + \tau_n^2 t^4/120] \]
and only use the $\sin$ and $\sinh$ formulas when $|\tau_n| h_n^2$ is 
sufficiently large.  

As a consequence, if we set
\[  y_n := \hat y(x_n), \quad y_n' := \hat y'(x_n) \]
for any $n$, then we have the forward recurrence
\begin{equation}    \label{FWRD}
  \left[ \begin{array}{c} y_{n+1} \\ y_{n+1}' \end{array} \right] = 
  \left[ \begin{array}{cc} \phi_n'(h_n) & \phi_n(h_n) \\ - \tau_n \phi_n(h_n) 
         & \phi_n'(h_n) \end{array} \right] 
  \left[ \begin{array}{c} y_n \\ y_n' \end{array} \right].
\end{equation}
If we denote the coefficient matrix in(\ref{FWRD}) by $A_n$, it is not 
difficult to show that it has inverse
\begin{equation}      \label{A-INV}
A_n^{-1} = \left[ \begin{array}{cc} \phi_n'(h_n) & -\phi_n(h_n) \\  
           \tau_n \phi_n(h_n) & \phi_n'(h_n) \end{array} \right].
\end{equation}
Hence, a backward recurrence is
\begin{equation}    \label{BWRD}
  \left[ \begin{array}{c} y_n \\ y_n' \end{array} \right] = 
  \left[ \begin{array}{cc} \phi_n'(h_n) & -\phi_n(h_n) \\ \tau_n \phi_n(h_n) 
         & \phi_n'(h_n) \end{array} \right] 
  \left[ \begin{array}{c} y_{n+1} \\ y_{n+1}' \end{array} \right].
\end{equation}
It can be seen when $\tau_n > 0$ that $A_n$ has eigenvalues
$\cos( \omega_n h_n) \pm i \sin( \omega_n h_n)$ and spectral radius one.
When $\tau_n < 0$, its eigenvalues are 
$\cosh( \omega_n h_n) \pm \sinh( \omega_n h_n)$ and its spectral radius is 
$\exp (\omega_n h_n)$.  The exponential factor reflects the potential 
mathematical instability of the initial value problem (\ref{SL_EQ}) 
when $\lambda < q(x)$.  To overcome this, define
\begin{equation}   \label{SIGMA}
 \sigma_n = \left\{ \begin{array}{ll} \exp(\omega_n h_n) & 
            \tau_n < -\epsilon \\  1 & {\rm otherwise},
            \end{array} \right.  
\end{equation}
and for $j \le k$
\begin{equation}   \label{PROD}
  p(j,k) = \sigma_j \sigma_{j+1} \cdots \sigma_k.
\end{equation}
Introduce the scaled variables
\begin{eqnarray}
 \tilde{y}_n &=& y_n / p(1, n-1)
                 \label{SCL1} \\
 \tilde{y}_n' &=& y_n' / p(1, n-1), 
                 \label{SCL2}
\end{eqnarray}
which satisfy the recurrences of the form (\ref{FWRD}) or (\ref{BWRD}) with 
coefficient matrix divided by $\sigma_n$.  These scaled matrices have spectral 
radius one.

To use (\ref{FWRD}) requires an initial condition to start, while 
(\ref{BWRD}) requires a terminal condition.  More generally, we define $u^F$
and $v^F$ to each satisfy the differential equation (\ref{yHat_EQ}) with 
respective initial conditions $u^F(0) = 1$, ${u^F}'(0) = 0$ and 
$v^F(0) = 0$, ${v^F}'(0) = 1$.   Similarly, $u^B$ and $v^B$ satisfy the same
differential equation but with respective terminal conditions  $u^B(\ell) = 1$, 
${u^B}'(\ell) = 0$ and $v^B(\ell) = 0$, ${v^B}'(\ell) = 1$.  We will define
the 2-vector $U^F$ to have components $u^F$ and ${u^F}'$; furthermore, for 
$n = 1, 2, \ldots, N+1$ let $U_n^F$ denote the two-vector with components 
$u_n^F$ and ${u_n^F}'$.  Define 2-vectors $U^B$, $V^F$, $V^B$, $U_n^B$, 
$V_n^F$, and $V_n^B$ analogously.  For the vectors with $B$ superscripts, we 
recur backwards from $n =N+1$ using (\ref{BWRD}) while for those with $F$
superscripts we recur forwards from $n=1$ using (\ref{FWRD}).  Finally, we 
use a tilde overscore ($\, \tilde{} \, $) to denote the scaled versions of 
these recurrences.   The analog of (\ref{SCL1})--(\ref{SCL2}) for $Y$ either 
$U$ or $V$ is
\begin{eqnarray}   
 \tilde{Y}_n^F &=& Y_n^F / p(1, n-1)
     \label{SCL3} \\
 \tilde{Y}_n^B &=& Y_n^B / p(n, N).  
     \label{SCL4}
\end{eqnarray}

Since $U^B$ and $V^B$ form a basis of solutions for (\ref{yHat_EQ}) there
exist constants $c_{11}, c_{12}, c_{21}, c_{22}$ such that
\begin{eqnarray}   
  U^F &=& c_{11} U^B + c_{12} V^B   \label{LC1} \\
  V^F &=& c_{21} U^B + c_{22} V^B   \label{LC2} 
\end{eqnarray}
for every $x$.  Define 
\begin{equation}    \label{DELTA}
   \Delta = u^B (v^B)' - (u^B)' v^B,
\end{equation}
then after some calculation it follows that
\begin{eqnarray}
  c_{11} &=& (u^F (v^B)' -  (u^F)' v^B) / \Delta  \label{c11} \\
  c_{12} &=& (u^B (u^F)' - (u^B)' u^F)  / \Delta  \label{c12} \\
  c_{21} &=& ((v^B)' v^F -  v^B (v^F)') / \Delta  \label{c21} \\
  c_{22} &=& (u^B (v^F)' - (u^B)' v^F)  / \Delta  \label{c22}
\end{eqnarray}
for any choice of $x \in [0, \ell]$.  

From the first component of (\ref{LC1})
\[  \hat u(\ell) = c_{11} u^B(\ell) + c_{12} v^B(\ell) = c_{11}.  \]
Similarly, 
\begin{eqnarray*}
  \hat u'(\ell) &=& c_{11} {u^B}'(\ell) + c_{12} {v^B}'(\ell) = c_{12}  \\
  \hat v(\ell)  &=& c_{21} u^B(\ell) + c_{22} v^B(\ell) = c_{21}  \\
  \hat v'(\ell) &=& c_{21} {u^B}'(\ell) + c_{22} {v^B}'(\ell) = c_{22} 
\end{eqnarray*}

Consequently, it follows from (\ref{fHat}) that
\begin{eqnarray}
  \hat f(\lambda) &=& \sfrac{\sqrt{\max\{0, 4 - [\hat u(\ell) + 
      \hat v'(\ell)]^2\}}} 
      {2 \pi |[\hat u'(\ell) \sin^2 \alpha
        + (\hat u(\ell) - \hat v'(\ell)) \sin \alpha \cos \alpha - 
           \hat v(\ell) \cos^2 \alpha]| }  \nonumber   \\
   &=& \sfrac{\sqrt{\max\{0, 4 - [c_{11} + c_{22}]^2\}}} 
       {2 \pi |c_{12} \sin^2 \alpha
        + (c_{11} - c_{22}) \sin \alpha \cos \alpha - 
           c_{21} \cos^2 \alpha|}.          \label{fHat1}
\end{eqnarray}
 
Since the formulas (\ref{c11})--(\ref{c22}) are valid for any $x$, 
another approach is to recur from both ends, computing the various 
$u_n^F$, $u_n^B$, $v_n^F$, and $v_n^B$, or their scaled equivalents, 
and `match' at some interior point denoted by $x = x_M$, to be determined 
below.  In detail, we begin with
\[ \tilde{u}_1^F = 1, \quad \tilde{v}_1^F = 0, \quad ({\tilde u}_1^F)' = 0, 
   \quad  ({\tilde v}_1^F)' = 1,  \]
\[ \tilde{u}_{N+1}^B = 1, \quad \tilde{v}_{N+1}^B = 0, \quad 
  ({\tilde u}_{N+1}^B)' = 0,  \quad  ({\tilde v}_{N+1}^B)' = 1,  \]
and then compute
\begin{eqnarray}
  \tilde{U}_{n+1}^F &=& A_n \tilde{U}_n^F / \sigma_n     \label{FWRDu}  \\
  \tilde{V}_{n+1}^F &=& A_n \tilde{V}_n^F / \sigma_n     \label{FWRDv}  
\end{eqnarray}
for $n = 1, 2, \ldots, M-1$, and
\begin{eqnarray}
  \tilde{U}_n^B &=& A_n^{-1} \tilde{U}_{n+1}^B / \sigma_n     \label{BWRDu} \\
  \tilde{V}_n^B &=& A_n^{-1} \tilde{V}_{n+1}^B / \sigma_n     \label{BWRDv} 
\end{eqnarray}
for $n = N, N-1, \ldots, M$.  Now from (\ref{SCL3})--(\ref{SCL4}) and 
(\ref{DELTA})--(\ref{c11}), with $n = M$ we have
\begin{eqnarray*}
  c_{11} &=& ((v_M^B)' u_M^F - v_M^B (u_M^F)') / (u_M^B (v_M^B)' 
               - (u_M^B)' v_M^B)   \\
  &=& p(1, M-1) (({\tilde{v}_M^B})' \tilde{u}_M^F - \tilde{v}_M^B 
                 ({\tilde{u}_M^F})') / [ p(M, N) (\tilde{u}_M^B 
                 ({\tilde{v}_M^B})' - ({\tilde{u}_M^B})' \tilde{v}_M^B)].
\end{eqnarray*}
Consequently, if we define the scale factor
\begin{equation}  \label{ZETA}
  \zeta_M := \sfrac{p(1, M-1)}{p(M,N)},
\end{equation}
then $c_{11}$ and similarly $c_{12}$, $c_{21}$, and $c_{22}$ given 
in (\ref{c11})--(\ref{c22}) must be multiplied by $\zeta_M$ if scaled 
variables are used.   To avoid rapid error buildup, it is desirable to have 
$\zeta_M \approx 1$, equivalently, 
\begin{eqnarray}   
  p(1, M-1) &\approx& p(M, N)    \nonumber \\
    &\approx& \sqrt{ p(1, N)}  \label{RATIO1}
\end{eqnarray}
with $M$ chosen to be an index for which the approximation is best.
In extreme cases, products of the $\sigma_n$ may overflow, so it is
best to work with their logs, i.e., from (\ref{SIGMA}), $h_n \omega_n$.
Then (\ref{RATIO1}) becomes
\begin{equation}    \label{RATIO2}
  \sum_{n=1}^{M-1}{}' \; h_n \omega_n \approx 0.5 \sum_{n=1}^N{}' \;
                      h_n \omega_n,
\end{equation}
where the $'$ on the sum means to replace $h_n \omega_n$ with zero for
any index corresponding to $\tau_n > 0$.   Care must also be taken in
scaling $c_{11}$, $c_{12}$, $c_{21}$, and $c_{22}$ by $\zeta$ to do 
the quotient with the logs first and only perform the exponentiation 
at the end.

Hence, for a given choice of $\lambda$, the stabilized algorithm first
makes an initial pass across the subintervals of $[0, \ell]$ to compute 
the $\{\sigma_n\}$ and $M$.  Next, the scaled forward and backward 
recurrences are performed that allow the computation of $c_{11}$, 
$c_{12}$, $c_{21}$, and $c_{22}$.  Finally (\ref{fHat1}) can be used 
to compute the estimates for $f(\lambda)$.  This can be repeated for a 
sequence of ever finer meshes until convergence is observed.  If this 
is not accomplished in a certain number of steps, the computation is 
suspended and an error flag is set.

As an illustration we choose the Mathieu equation for which
\begin{equation}    \label{MATHIEU}
   q(x) = \cos x.  
\end{equation}
We computed the spectral function $f(\lambda)$ at 101 equally spaced
$\lambda$ values in several stability intervals, with simple and 
with double shooting; the average time was measured for each method and 
interval.  In all cases a hundred repetitions were made for each of the 
$\lambda$ values in order for the computer clock to produce a reliable 
time.  This was done at several tolerances on the $\hat f$ sequence.
The output is summarized in Table~3.1 for a Dirichlet condition 
at $x = 0$ corresponding to the choice $\alpha = 0$.  Two absolute error 
tolerances were used: $10^{-6}$ and $10^{-8}$.    The `failure' column 
shows a count of the number of values for which convergence was not 
achieved in eight mesh refinements (bisected uniform meshes).  In 
Table~3.2 are the corresponding values for a Neumann condition 
($\alpha = \pi / 2$). 

\vskip 18pt

{\small
\begin{center}
Table 4.1. Simple vs. double shooting for (\ref{MATHIEU}) -- Dirichlet.

\vskip 8pt
\begin{tabular}{ccccccc}
\hline

&&\multicolumn{2}{c}{simple shooting}& &\multicolumn{2}{c}{double shooting}\\
tolerance&interval&  time   & \# failures& { } &  time  & \# failures \\
\hline
\\
$10^{-6}$&$[-0.3784,-0.3476]$& 0.253 & 42 & { } & 0.023 &  0 \\
&$[ 0.5949, 0.9180]$& 0.080 &  0 &     & 0.007 &  0 \\
&$[ 1.2932, 2.2851]$& 0.006 &  0 &     & 0.007 &  0 \\
&$[ 2.3426, 4.0319]$& 0.006 &  0 &     & 0.006 &  0 \\
\\
$10^{-8}$&$[-0.3784,-0.3476]$& 0.659 & 96 & { } & 0.042 &  0 \\
&$[ 0.5949, 0.9180]$& 0.458 & 51 &     & 0.013 &  0 \\
&$[ 1.2932, 2.2851]$& 0.021 &  0 &     & 0.021 &  0 \\
&$[ 2.3426, 4.0319]$& 0.017 &  0 &     & 0.017 &  0 \\
\end{tabular}  
\end{center}
} 
\vskip 18pt


{\small
\begin{center}
Table 4.2. Simple vs. double shooting for (\ref{MATHIEU}) -- Neumann.

\vskip 8pt
\begin{tabular}{ccccccc}
\hline

&&\multicolumn{2}{c}{simple shooting}& &\multicolumn{2}{c}{double shooting}\\
tolerance&interval&  time   & \# failures& { } &  time  & \# failures \\
\hline
\\
$10^{-6}$&$[-0.3784,-0.3476]$& 0.300 & 47 & { } & 0.023 &  0 \\
&$[ 0.5949, 0.9180]$& 0.110 &  1 &     & 0.007 &  0 \\
&$[ 1.2932, 2.2851]$& 0.007 &  0 &     & 0.007 &  0 \\
&$[ 2.3426, 4.0319]$& 0.009 &  0 &     & 0.009 &  0 \\
\\
$10^{-8}$&$[-0.3784,-0.3476]$& 0.590 & 96 & { } & 0.040 &  0 \\
&$[ 0.5949, 0.9180]$& 0.491 & 68 &     & 0.016 &  0 \\
&$[ 1.2932, 2.2851]$& 0.018 &  0 &     & 0.018 &  0 \\
&$[ 2.3426, 4.0319]$& 0.021 &  0 &     & 0.021 &  0 \\
\end{tabular}  
\end{center}
} 
\vskip 18pt

Clearly the simple shooting approach falters on the first two
stability intervals.  For the other two intervals the quantities 
$\tau_n$ in (\ref{TAU}) are always positive so that theoretically 
the two methods should be equally reliable.  The output supports 
this, and little time is lost from the minor overhead of the 
double shooting.  For the remainder of the output in this paper 
the double shooting method always will be used.

\section{Indeterminate cases}                        
\setcounter{equation}{0}

When the potential $q$ is periodic, it is known that the spectrum exhibits
spectral gaps of resolvent set where no spectrum can occur, i.e, where
$f(\lambda) = 0$.  Moreover, the endpoints of a spectral gap occur at
values of $\lambda^*$ for which $(u + v')(\ell, \lambda^*) = \pm 2$, i.e.,
the numerator in (\ref{f4}) vanishes.  For some examples the denominator 
may also vanish at the same $\lambda^*$.   In such cases we might expect 
the numerical error to be large for values of $\lambda$ near such 
$\lambda^*$.  In fact, such $\lambda^*$ arise at endpoints of spectral 
gaps for any potential exhibiting even symmetry, i.e., $q(\ell - x) = q(x)$ 
for all $x$, as is the case for Mathieu's equation.

\vskip 8pt \noindent
{\bf Case 1 (Dirichlet)}: assume that for some fixed $\lambda = \lambda^*$ 
we have $v(\ell, \lambda^*) = 0$ and $u(\ell, \lambda^*) + v_x(\ell, 
\lambda^*) = \pm 2$.   Then, near $\lambda^*$ we have
\begin{eqnarray*}
  4 - [u(\ell, \lambda) + v_x(\ell, \lambda)]^2 &=& [2 + |u(\ell, \lambda) 
   + v_x(\ell, \lambda)|][2 - |u(\ell, \lambda) + v_x(\ell, \lambda)|] \\
   &=& [2 + |u(\ell, \lambda) + v_x(\ell, \lambda)|]
       [(u_{\lambda}(\ell, \lambda^*) + v_{x \lambda}(\ell, \lambda^*)) 
       (\lambda - \lambda^*) + O((\lambda - \lambda^*)^2)] 
\end{eqnarray*}
and 
\[  v(\ell, \lambda) = v_{\lambda}(\ell, \lambda^*) (\lambda - \lambda^*) 
    + O((\lambda- \lambda^*)^2),  \]
so that for $\lambda < \lambda^*$
\begin{eqnarray}  
   f(\lambda) &=& \sfrac{\sqrt{[2 + |u(\ell, \lambda) + v_x(\ell, \lambda)|]
     |u_{\lambda}(\ell, \lambda^*) + v_{x \lambda}(\ell, \lambda^*)|}} 
      {2\pi |v_{\lambda}(\ell, \lambda^*)| \sqrt{\lambda^* - \lambda}} 
      + O(\lambda^* - \lambda)    \nonumber \\
   &\approx& \sfrac{\sqrt{[2 + |u(\ell, \lambda) + v_x(\ell, \lambda)|]
     |u_{\lambda}(\ell, \lambda^*) + v_{x \lambda}(\ell, \lambda^*)|}} 
      {2\pi |v_{\lambda}(\ell, \lambda^*)| \sqrt{\lambda^* - \lambda}}.    
      \label{FIX_D}
\end{eqnarray} 

\vskip 8pt \noindent
{\bf Case 2: Neumann}: assume that for some fixed $\lambda = \lambda^*$ 
we have $u_x(\ell, \lambda^*) = 0$ and $u(\ell, \lambda^*) + 
v_x(\ell, \lambda^*) = \pm 2$.  Analogous to the Dirichlet case, we have 
for $\lambda > \lambda^*$ 
\begin{eqnarray}   
   f(\lambda) &=& \sfrac{\sqrt{[2 + |u(\ell, \lambda) + v_x(\ell, \lambda)|]
      |u_{\lambda}(\ell, \lambda^*) + v_{x \lambda} (\ell, \lambda^*)|}} 
      {2\pi |u_{x \lambda}(\ell, \lambda^*)| \sqrt{\lambda - \lambda^*}} 
      + O(\lambda- \lambda^*)    \nonumber \\
      &\approx& \sfrac{\sqrt{[2 + |u(\ell, \lambda) + v_x(\ell, \lambda)|]
      |u_{\lambda}(\ell, \lambda^*) + v_{x \lambda} (\ell, \lambda^*)|}} 
      {2\pi |u_{x \lambda}(\ell, \lambda^*)| \sqrt{\lambda - \lambda^*}}.
      \label{FIX_N}
\end{eqnarray} 

For a step-function potential the partial derivatives appearing in the above
$f$ formulas can be computed easily from the closed form solutions given
in the previous section.    The details are given in an appendix.
Note that in either case we expect a $1 / \sqrt{|\lambda - \lambda^*|}$
behavior near a point of indeterminacy $\lambda^*$.   Our experience has 
shown that the expected loss of significance is not serious except (1) at 
very tight tolerances, (2) at values of $\lambda$ very close to 
$\lambda^*$, or (3) near gap endpoints where the gap is very narrow 
(larger $\lambda^*$).   

As an illustration we again choose the Mathieu potential (\ref{MATHIEU}).  
For an absolute error tolerance of $10^{-8}$, we evaluated (\ref{f4})
and (\ref{FIX_D}) near endpoints of the stability intervals.  For
(\ref{f4}) we also estimated the rate $\alpha$ in 
\[ f(\lambda) \approx \sfrac{\mbox{constant}}{|\lambda - \lambda^*|^{\alpha}} \]
by
\[ \mbox{rate} \approx \sfrac{\log [f(\lambda_2) / f(\lambda_1)]}
   {\log[|\lambda_1 - \lambda^*|/|\lambda_2 - \lambda^*|]} .  \]
Table~5.1 displays the numerical output for a Dirichlet initial condition, 
where the indeterminacy occurs at the right-hand end of a stability interval.
The respective $\lambda^*$ values for (\ref{FIX_D}) are 
\[ \{-0.347669125306, 0.918058176625, 2.28515693444\}. \]
Table~5.2 does the same for a Neumann initial condition, where the indeterminacy 
is at the left-hand end of a stability interval.  The $\lambda^*$ values are 
\[  \{-0.378489221265, 0.594799970122, 1.29316628334\}.  \]

\vskip 18pt
\pagebreak
{\small
\begin{center}
Table 5.1. Behavior near an indeterminacy for (\ref{MATHIEU}) -- Dirichlet.

\vskip 8pt
\begin{tabular}{cccr}
\hline
$\lambda$ & $f$ from (\ref{f4}) & rate & $f$ from (\ref{FIX_D}) \\
\hline
$-$0.3497 & 1.34079 &       & 1.38601 \\  
$-$0.3493 & 1.50630 & 0.531 & 1.54665 \\  
$-$0.3489 & 1.74540 & 0.524 & 1.78029 \\
$-$0.3485 & 2.13833 & 0.517 & 2.16680 \\
$-$0.3481 & 2.98860 & 0.510 & 3.00872 \\
$-$0.3477 &11.23586 & 0.514 &11.21862 \\
\\
   0.9157 & 2.35819 &       & 2.35955 \\  
   0.9161 & 2.58811 & 0.500 & 2.58938 \\
   0.9165 & 2.90162 & 0.500 & 2.90281 \\
   0.9169 & 3.36590 & 0.500 & 3.36703 \\
   0.9173 & 4.16048 & 0.500 & 4.16167 \\
   0.9177 & 6.05367 & 0.500 & 6.05564 \\
\\
   2.2831 & 1.90694 &       & 1.87531 \\
   2.2835 & 2.11789 & 0.485 & 2.08944 \\
   2.2839 & 2.42382 & 0.488 & 2.39896 \\
   2.2843 & 2.92599 & 0.492 & 2.90536 \\
   2.2847 & 3.99392 & 0.495 & 3.97862 \\
   2.2851 &11.27750 & 0.498 &11.26558
\end{tabular}  
\end{center}
} 
\vskip 18pt

{\small
\begin{center}
Table 5.2. Behavior near an indeterminacy for (\ref{MATHIEU}) -- Neumann.

\vskip 8pt
\begin{tabular}{cccr}
\hline
$\lambda$ & $f$ from (\ref{f4}) & rate & $f$ from (\ref{FIX_N}) \\
\hline
$-$0.3784 & 5.21621 & 0.504 & 5.22961 \\
$-$0.3780 & 2.21324 & 0.511 & 2.23105 \\
$-$0.3476 & 1.63105 & 0.518 & 1.65468 \\
$-$0.3472 & 1.34574 & 0.525 & 1.37416 \\
$-$0.3468 & 1.16787 & 0.532 & 1.20046 \\
$-$0.3464 & 1.04307 &       & 1.07943 \\
\\
   0.5952 &10.46971 & 0.500 &10.47355 \\
   0.5956 & 7.40089 & 0.501 & 7.40592 \\
   0.5960 & 6.04084 & 0.501 & 6.04691 \\
   0.5964 & 5.22980 & 0.501 & 5.23678 \\
   0.5968 & 4.67613 & 0.502 & 4.68392 \\
   0.5972 & 4.26729 &       & 4.27581 \\
\\
   1.2936 &10.78605 & 0.500 &10.77975 \\
   1.2940 & 7.78143 & 0.500 & 7.77625 \\
   1.2944 & 6.39829 & 0.499 & 6.39287 \\
   1.2948 & 5.56142 & 0.499 & 5.55555 \\
   1.2952 & 4.98576 & 0.499 & 4.97941 \\
   1.2956 & 4.55873 &       & 4.55190 \\
\end{tabular}  
\end{center}
}

\pagebreak 
\vskip 18pt

There are slight differences between the two approaches.  Since           
it is difficult to calculate exact answers in these cases (and $\lambda^*$ 
itself), we have no easy way to judge which, if either, is more correct.
By an inspection of intermediate quantities needed for the special
formula (\ref{FIX_D}), viz., $u_{\lambda}$, $v_{x,\lambda}$, they can
be quite sensitive to the error in $|\lambda^* - \lambda|$, as well as
the tolerance. In a more positive vein, it is clear that the double shooting 
method is in agreement as to the growth rate of $f$ near $\lambda^*$ in 
the indeterminate situations.

\section{Other Numerical results}                  
\setcounter{equation}{0}

In this section we exhibit computational results illustrating the algorithms
developed in the previous sections.  For brevity we choose five potentials;
the first is Mathieu's equation (\ref{MATHIEU}).  As mentioned in the previous 
section, potentials such as this one that exhibit even symmetry have special 
properties.  It can be shown that if 
$q(\ell - x) = q(x)$ for every $x$, then
\begin{equation}   \label{EVEN1}
   u(\ell, \lambda) = v'(\ell, \lambda) \qquad \mbox{for every } \lambda.
\end{equation}
Moreover, whenever $|u(\ell, \lambda^*)| = 1$ for some $\lambda^*$, then
either 
\begin{equation}   \label{EVEN2}
   v(\ell, \lambda^*) = 0 
\end{equation}
or
\begin{equation}   \label{EVEN3}
   u'(\ell, \lambda^*) = 0.
\end{equation}
In the case (\ref{EVEN2}), $\lambda^*$ is the left endpoint of a spectral gap 
when $y$ in (\ref{SL_EQ}) satisfies a Dirichlet condition.  In the case of
(\ref{EVEN3}),  $\lambda^*$ is the right endpoint of a spectral gap 
when $y$ in (\ref{SL_EQ}) satisfies a Neumann condition.

The potentials in our other examples are 
\begin{eqnarray}   
    & & 3 / (2 + \sin x),              \label{EX2}    \\
    & & 1 / \sqrt{1 -0.75 \sin^2 x},   \label{EX3}   \\
    & & (0.5 + \cos x + \cos 2x + \cos 3x) / \pi,   \label{EX4}  \\
    & & \sin x + 0.5 \sin 2x + 0.1 \sin 3x.   \label{EX5}
\end{eqnarray}
These have period $\ell = 2 \pi$ except for (\ref{EX3}) that has period
$\pi$.  In addition to (\ref{MATHIEU}), examples (\ref{EX3}) and (\ref{EX4}) 
also have even symmetry.

In Table~6.1a we display the endpoints of the first few stability intervals 
for the first two examples.  For Mathieu's equation these are known 
\cite{INCE} from the theory of elliptic cylinder functions.  The numerical 
values agree with those found in \cite{INCE}, or see \cite[Table I]{TOMSEAST}.  
The numerical method used was simple binary search (bisection) seeking the 
zeros of
\begin{equation}    \label{GMU}
  g(\lambda) := 2 - |u(\ell, \lambda) + v'(\ell, \lambda)|.
\end{equation}
An absolute error tolerance of $10^{-8}$ was used in all cases.  Values of 
$f$ were first computed over a sufficiently fine grid to identify the 
locations of the gaps.  As $\lambda$ increases the gap width narrows, making 
it more difficult to isolate gap boundaries.  Moreover, the loss of 
significance in evaluating $g$ worsens; eventually we may have to switch to
the techniques in Section~4 to help overcome this.    However, this was not 
necessary for the data in Table~5.1a.

\vskip 18pt
{\small
\begin{center}
Table 6.1a. Stability intervals for the first two examples.

\vskip 8pt
\begin{tabular}{rrcrr}
\hline
\multicolumn{2}{c}{Mathieu}& & \multicolumn{2}{c}{Example 5.4} \\
\hline
($-$0.378489, & $-$0.347669) & { } &( 2.250000, & 2.548882) \\
(   0.594800, &    0.918058) &     &( 3.055360, & 3.941647) \\
(   1.293166, &    2.285157) &     &( 4.146186, & 5.736211) \\
(   2.342581, &    4.031922) &     &( 5.796032, & 7.994726) \\
(   4.035301, &    6.270837) &     &( 8.010349, &10.743819) \\
(   6.270945, &    9.014297) &     &(10.747778, &13.991464) \\
\hline
\end{tabular}  
\end{center}
} 
\vskip 18pt

Similarly, Table~6.1b contains the stability intervals for Examples 
(\ref{EX3})--(\ref{EX5}).  

\vskip 18pt
{\small
\begin{center}
Table 6.1b. Stability intervals for the last three examples.

\vskip 8pt
\begin{tabular}{rrcrrcrr}
\hline
\multicolumn{2}{c}{Example 5.5}& & \multicolumn{2}{c}{Example 5.6} 
& & \multicolumn{2}{c}{Example 5.7} \\
\hline
( 1.346160, &    2.136962) & { } & (0.106301, & 0.247914) & { } & ($-$0.419549, & $-$0.391618) \\
( 2.594046, &    5.310602) &     & (0.503181, & 0.995282) &     & (   0.570873, &    0.840333) \\
( 5.452072, &   10.356984) &     & (1.311604, & 2.240365) &     & (   1.362407, &    2.217768) \\
(10.396276, &   17.369252) &     & (2.602473, & 4.151030) &     & (   2.442559, &    4.011052) \\
(17.380456, &   26.372454) &     & (4.198967, & 6.407883) &     & (   4.078880, &    6.271355) \\
(26.375745, &   37.373218) &     & (6.426576, & 9.160844) &     & (   6.283327, &    9.017477) \\
\hline
\end{tabular}  
\end{center}
} 
\vskip 18pt

Next we compare the new formula (\ref{f4}) with a variant of the SLEDGE 
code.  The original SLEDGE \cite{TOMSEAST}, \cite{TOMS}, \cite{SLEDGE} 
could return estimates for the spectral measure $\rho(\lambda)$
but not for the density function $f(\lambda)$; to this code was added an 
implementation of interpolant~3 from \cite[Eqn. (3.5)]{INTERP}, there 
denoted by $(I_3 \rho_b)'$, in order to provide estimates for $f(\lambda)$. 
Recall that SLEDGE uses the Levitan-Levinson characterization of the measure 
$\rho(\lambda)$.  This requires the calculation of many eigenvalues and suitably 
normalized eigenfunctions of a regularized Sturm-Liouville problem over a finite 
interval $(0, b)$ for a sequence of increasingly larger $b$.  Nevertheless,
it was demonstrated in \cite{TOMSEAST} that SLEDGE is capable of successfully
handling a wide scope of problems.  Our first example for Table~6.2a is the 
Mathieu equation (\ref{MATHIEU}) with a Dirichlet initial condition
($\alpha = 0$), in Table~6.2b are the data corresponding to $\alpha = \pi /6$, 
and in Table~6.2c are the data corresponding to a Neumann initial condition.   
The internal tolerance used by SLEDGE in the calculation of the eigenvalues and 
eigenfunctions was $10^{-6}$, while for the new method it was either $10^{-6}$
or $10^{-8}$, as shown. The final line of the table shows the computer time 
needed for computing $f$ and $\rho$ at 66 $\lambda$ points.  Since SLEDGE has 
no choice but to compute both $\rho$ and $f$, we required our new method do 
both as well.  For brevity, only some $f$ output values and no $\rho$ values 
are shown in the tables.

SLEDGE loses accuracy near the boundaries of the stability intervals,
though it still seems to be converging as $b$ increases.  There is 
little difference in using the new formula (\ref{f4}) at the tighter 
tolerance, other than an increase in time.  What differences there 
are generally occur near endpoints of stability intervals, especially 
when the formulas are indeterminate there.  In all cases it is clear 
that the new approach is much faster.

\pagebreak
\vskip 0.2in
\begin{center}
Table 6.2a. Estimates of $f(\lambda)$ for Mathieu's Equation -- Dirichlet.
\vskip 4pt

\begin{tabular}{rcccc}
\hline
$\lambda  $&SLEDGE&SLEDGE& (\ref{f4}) & (\ref{f4}) \\
\hline
$-$0.38 & 0.000000 & 0.000000 & 0.00000000 & 0.00000000 \\
$-$0.37 & 0.221836 & 0.221583 & 0.22149618 & 0.22149622 \\
$-$0.36 & 0.439526 & 0.438250 & 0.43801179 & 0.43801181 \\
$-$0.35 & 1.262850 & 1.247271 & 1.24515689 & 1.24515701 \\
$-$0.34 & 0.001181 & 0.000591 & 0.00000000 & 0.00000000 \\
\multicolumn{5}{c}{spectral gap} \\
  0.59  & 0.000029 & 0.000006 & 0.00000000 & 0.00000000 \\
  0.60  & 0.034764 & 0.034936 & 0.03503180 & 0.03503178 \\
  0.70  & 0.176067 & 0.175985 & 0.17595739 & 0.17595738 \\
  0.80  & 0.305406 & 0.304771 & 0.30451657 & 0.30451657 \\
  0.90  & 0.876650 & 0.855793 & 0.84810870 & 0.84810870 \\
  0.92  & 0.025007 & 0.012603 & 0.00000000 & 0.00000000 \\
\multicolumn{5}{c}{spectral gap} \\
  1.29  & 0.000896 & 0.000195 & 0.00000000 & 0.00000000 \\
  1.30  & 0.028250 & 0.036852 & 0.03714803 & 0.03714802 \\
  1.50  & 0.188818 & 0.188736 & 0.18871550 & 0.18871550 \\
  1.75  & 0.268937 & 0.268931 & 0.26892936 & 0.26892936 \\
  2.00  & 0.336884 & 0.336798 & 0.33675478 & 0.33675478 \\
  2.25  & 0.673109 & 0.591205 & 0.56713172 & 0.56713172 \\
  2.28  & 0.391315 & 0.839088 & 1.23367035 & 1.23367034 \\
  2.30  & 0.267555 & 0.179038 & 0.00000000 & 0.00000000 \\
\multicolumn{5}{c}{spectral gap} \\
  2.34  & 0.184005 & 0.071690 & 0.00000000 & 0.00000000 \\
  2.50  & 0.330775 & 0.330639 & 0.33062792 & 0.33062792 \\
  2.75  & 0.392572 & 0.392579 & 0.39258965 & 0.39258966 \\
  3.00  & 0.431331 & 0.431330 & 0.43133033 & 0.43133034 \\
  3.25  & 0.463740 & 0.463753 & 0.46375540 & 0.46375540 \\
  3.50  & 0.493118 & 0.493118 & 0.49311812 & 0.49311813 \\
\hline
  $b$     &$64 \pi$& $128 \pi$&            &            \\ 
tolerance &$10^{-6}$&$10^{-6}$&  $10^{-6}$ &  $10^{-8}$ \\
total time&  72.67 &  246.35  &    0.10    &    0.17   
\end{tabular}
\end{center}
\vskip 0.2in

\pagebreak
\vskip 0.2in
\begin{center}
Table 6.2b. Estimates of $f(\lambda)$ for Mathieu's Equation -- $\alpha = \pi / 6$.
\vskip 4pt

\begin{tabular}{rcccc}
\hline
$\lambda  $&SLEDGE&SLEDGE& (\ref{f4}) & (\ref{f4}) \\
\hline
$-$0.38 & 0.000000 & 0.000000 & 0.00000000 & 0.00000000 \\
$-$0.37 & 0.254422 & 0.254359 & 0.25428685 & 0.25428585 \\
$-$0.36 & 0.358362 & 0.358111 & 0.35803367 & 0.35803367 \\
$-$0.35 & 0.271701 & 0.271974 & 0.27213584 & 0.27213584 \\
$-$0.34 & 0.000000 & 0.177284 & 0.00000000 & 0.00000000 \\
\multicolumn{5}{c}{spectral gap} \\
  0.59  & 0.000000 & 0.001934 & 0.00000000 & 0.00000000 \\
  0.60  & 0.046178 & 0.046391 & 0.04652124 & 0.04652122 \\
  0.70  & 0.213029 & 0.212952 & 0.21292212 & 0.21292210 \\
  0.80  & 0.311584 & 0.311230 & 0.31111121 & 0.31111121 \\
  0.90  & 0.334158 & 0.336703 & 0.33591487 & 0.33591485 \\
  0.92  & 0.064154 & 0.065358 & 0.00000000 & 0.00000000 \\
\multicolumn{5}{c}{spectral gap} \\
  1.29  & 0.018331 & 0.015724 & 0.00000000 & 0.00000000 \\
  1.30  & 0.036250 & 0.048914 & 0.04930685 & 0.04930685 \\
  1.50  & 0.225283 & 0.225239 & 0.22523166 & 0.22523166 \\
  1.75  & 0.289651 & 0.289648 & 0.28965416 & 0.28965416 \\
  2.00  & 0.327035 & 0.327015 & 0.32700588 & 0.32700588 \\
  2.25  & 0.365013 & 0.371550 & 0.36740588 & 0.36740588 \\
  2.28  & 0.550911 & 0.173426 & 0.27382995 & 0.27382995 \\
  2.30  & 1.003042 & 0.719051 & 0.00000000 & 0.00000000 \\
\multicolumn{5}{c}{spectral gap} \\
  2.34  & 0.604674 & 0.435577 & 0.00000000 & 0.00000000 \\
  2.50  & 0.324044 & 0.324191 & 0.32423291 & 0.32423291 \\
  2.75  & 0.347326 & 0.347336 & 0.34733465 & 0.34733465 \\
  3.00  & 0.356739 & 0.356748 & 0.35675152 & 0.35675152 \\
  3.25  & 0.362111 & 0.362122 & 0.36212172 & 0.36212172 \\
  3.50  & 0.365271 & 0.365276 & 0.36527626 & 0.36527626 \\
\hline
  $b$     &$64 \pi$& $128 \pi$&            &            \\ 
tolerance &$10^{-6}$&$10^{-6}$&  $10^{-6}$ &  $10^{-8}$ \\
total time&  81.45 &  231.55  &    0.09    &    0.14   
\end{tabular}
\end{center}
\vskip 0.2in

\pagebreak
\vskip 0.2in
\begin{center}
Table 6.2c. Estimates of $f(\lambda)$ for Mathieu's Equation -- Neumann.
\vskip 4pt

\begin{tabular}{rcccc}
\hline
$\lambda  $&SLEDGE&SLEDGE& (\ref{f4}) & (\ref{f4}) \\
\hline
$-$0.38 & 0.000000 & 0.000000 & 0.00000000 & 0.00000000 \\
$-$0.37 & 0.458625 & 0.457836 & 0.45743969 & 0.45743978 \\
$-$0.36 & 0.231738 & 0.231458 & 0.23132070 & 0.23132067 \\
$-$0.35 & 0.081480 & 0.081348 & 0.08137240 & 0.08137222 \\
$-$0.34 & 0.000162 & 0.000080 & 0.00000000 & 0.00000000 \\
\multicolumn{5}{c}{spectral gap} \\
  0.59  & 0.019597 & 0.009810 & 0.00000000 & 0.00000000 \\
  0.60  & 2.988328 & 2.893315 & 2.89226446 & 2.89226447 \\
  0.70  & 0.577118 & 0.576146 & 0.57582799 & 0.57582799 \\
  0.80  & 0.333326 & 0.332866 & 0.33272798 & 0.33272798 \\
  0.90  & 0.119832 & 0.119566 & 0.11946722 & 0.11946721 \\
  0.92  & 0.000890 & 0.003814 & 0.00000000 & 0.00000000 \\
\multicolumn{5}{c}{spectral gap} \\
  1.29  & 0.113623 & 0.057586 & 0.00000000 & 0.00000000 \\
  1.30  & 2.917653 & 2.769684 & 2.72749873 & 2.72749865 \\
  1.50  & 0.538198 & 0.537216 & 0.53689910 & 0.53689910 \\
  1.75  & 0.376904 & 0.376796 & 0.37675762 & 0.37675762 \\
  2.00  & 0.300894 & 0.300884 & 0.30087526 & 0.30087526 \\
  2.25  & 0.162234 & 0.180877 & 0.17865547 & 0.17865547 \\
  2.28  & 0.062850 & 0.040850 & 0.08212985 & 0.08212987 \\
  2.30  & 0.160618 & 0.097636 & 0.00000000 & 0.00000000 \\
\multicolumn{5}{c}{spectral gap} \\
  2.34  & 0.356154 & 0.299811 & 0.00000000 & 0.00000000 \\
  2.35  & 0.405038 & 0.440822 & 0.82126926 & 0.82126926 \\
  2.50  & 0.307667 & 0.306737 & 0.30645078 & 0.30645078 \\
  2.75  & 0.258179 & 0.258104 & 0.25808419 & 0.25808419 \\
  3.00  & 0.234940 & 0.234913 & 0.23490391 & 0.23490391 \\
  3.25  & 0.218489 & 0.218484 & 0.21847979 & 0.21847979 \\
  3.50  & 0.205484 & 0.205474 & 0.20547041 & 0.20547041 \\
\hline
  $b$     &$64 \pi$& $128 \pi$&            &            \\ 
tolerance &$10^{-6}$&$10^{-6}$&  $10^{-6}$ &  $10^{-8}$ \\
total time&  73.64 &  239.55  &    0.09    &    0.16   
\end{tabular}
\end{center}
\vskip 0.2in


\section{Appendix: estimating variational quantities}      
\setcounter{equation}{0}

Here we derive a method for computing the partial derivative with respect
to $\lambda$ of the quantities given in Section~5 for overcoming
indeterminacies.  Recall from Section~4 that the forward recurrence is
\begin{equation}   \label{U_FWRD}
  U_{n+1}^F = A_n U_n^F 
\end{equation}
so that
\begin{equation}   \label{U_BWRD}
  U_{n+1, \lambda}^F = A_{n,\lambda} U_n^F + A_n U_{n, \lambda}^F.
\end{equation}
Omitting the $n$ subscripts for now, we have 
\[  A_{\lambda} = \left[ \begin{array}{cc} \phi_{x \lambda} & 
        \phi_{\lambda} \\  -\tau_{\lambda} \phi - \tau 
        \phi_{\lambda} & \phi_{x \lambda} \end{array} \right].  \]
But from (\ref{TAU})
\[  \tau_{\lambda} = 1;  \]
furthermore, for either sign on $\tau$ it is easily shown that at 
$x = x_{n+1}$
\[  \phi_{x \lambda} = - h \phi / 2, \]
and
\[  \phi_{\lambda} = (h \phi_x  - \phi) / (2 \tau).  \]
Consequently, with the $n$ subscripts restored, a forward recursion is
\begin{eqnarray}   
  \sfrac{\partial U_{n+1}^F}{\partial \lambda} &=& 
  0.5 \left[ \begin{array}{cc} 
       -h_n \phi_n(h_n) & (h_n \phi_{n,x}(h_n) - \phi_n(h_n)) / \tau_n \\ 
       -\phi_n(h_n) - h_n \phi_{n, x}(h_n) & -h_n \phi_n(h_n) 
        \end{array} \right]  U_n^F  \nonumber \\
    & & + \left[ \begin{array}{cc} \phi_{n,x}(h_n) & \phi_n(h_n) \\ 
       -\tau_n \phi_n(h_n) & \phi_{n,x}(h_n) \end{array} \right] 
        \sfrac{\partial U_n^F}{\partial \lambda}  \label{VAR_FWRD_N}
\end{eqnarray}
with
\[  \sfrac{\partial U_1^F}{\partial \lambda} = \left[ \begin{array}{c} 0 
\\ 0 \end{array} \right].  \]
The forward recurrence for $V_{n, \lambda}^F$ is identical -- only the 
initial conditions on $V_1^F$ differ.

A similar analysis using the $A_n^{-1}$ as given in (\ref{A-INV})
leads to the backward recurrence
\begin{eqnarray}
  \sfrac{\partial U_n^B}{\partial \lambda} &=& 
  0.5 \left[ \begin{array}{cc} 
       -h_n \phi_n(h_n) & (\phi_n(h_n) - h_n \phi_{n,x}(h_n)) / \tau_n \\
        \phi_n(h_n) + h_n \phi_{n, x}(h_n) &  -h_n \phi_n(h_n) 
      \end{array} \right]  U_{n+1}^B  \nonumber \\
    & & + \left[ \begin{array}{cc} \phi_{n,x}(h_n) & -\phi_n(h_n) \\ 
       \tau_n \phi_n(h_n) & \phi_{n,x}(h_n) \end{array} \right] 
       \sfrac{\partial U_{n+1}^B}{\partial \lambda}  \label{VAR_BWRD_N}
\end{eqnarray}
with
\[  \sfrac{\partial U_{N+1}^B}{\partial \lambda} = \left[ \begin{array}{c} 0 
\\ 0 \end{array} \right].  \]
Again this holds as well for $V_{n, \lambda}^B$ with appropriate terminal values
for $V_{N+1}^B$.

As we saw in Section~4 it is desirable to scale the variables.  We will
use the same notation as \S4, and will first develop the formulas for 
$\tilde{U}^F$ as those for $\tilde{U}^B$, $\tilde{V}^F$, and 
$\tilde{V}^B$ are analogous.   Following (\ref{SCL3})--(\ref{SCL4}), we 
define
\begin{eqnarray*}   
 \tilde{U}_n^F &=& U_n^F / p(1, n-1)  \\
 \tilde{U}_n^B &=& U_n^B / p(n, N) .  
\end{eqnarray*}
Differentiation with respect to $\lambda$ yields
\[  \sfrac{\partial {\tilde U}_n^F}{\partial \lambda} = \sfrac{\partial U_n^F}
   {\partial \lambda} / p(1,n-1) - U_n^F \sfrac{\partial p(1,n-1)}
   {\partial \lambda} / p(1,n-1)^2. \]
But
\begin{eqnarray*}
   -\sfrac{\partial p(1,n-1)}{\partial \lambda} / p(1,n-1)^2 &=& [0.5 / p(1,n-1)^2] \sum \sfrac{h_j p(1,n-1)}
         {\sqrt{-\tau_j}}   \\
     &=& \sum (h_j / \sqrt{-\tau_j}) / (2 p(1,n-1)),
\end{eqnarray*}
where the sum is taken over $j$, $1 \le j < n$, for which $\tau_j < 0$.
Hence, 
\begin{equation}
    \sfrac{\partial {\tilde U}_n^F}{\partial \lambda} = \sfrac{1}{p(1,n-1)} 
    \left[ \sfrac{\partial U_n^F}{\partial \lambda} + U_n^F \sum \sfrac{h_j}
    {\sqrt{-\tau_j}} \right],
\end{equation}
and similarly
\begin{eqnarray}
    \sfrac{\partial {\tilde V}_n^F}{\partial \lambda} &=& \sfrac{1}{p(1,n-1)} 
    \left[ \sfrac{\partial V_n^F}{\partial \lambda} + V_n^F \sum \sfrac{h_j}
    {\sqrt{-\tau_j}} \right] \\
    \sfrac{\partial {\tilde U}_n^B}{\partial \lambda} &=& \sfrac{1}{p(n,N)} 
    \left[ \sfrac{\partial U_n^B}{\partial \lambda} + U_n^B \sum \sfrac{h_j}
    {\sqrt{-\tau_j}} \right] \\
    \sfrac{\partial {\tilde V}_n^B}{\partial \lambda} &=& \sfrac{1}{p(n,N)} 
    \left[ \sfrac{\partial V_n^B}{\partial \lambda} + V_n^B \sum \sfrac{h_j}
    {\sqrt{-\tau_j}} \right]. \\
\end{eqnarray}
The sums for the $F$-superscripted cases run from $j=1$ to $j=n-1$, while 
those for the $B$-superscripted cases go from $j=n$ to $j=N$.  In either
case, the indices for which $\tau_j > 0$ are omitted.  Recall that we expect 
a more stable algorithm if we recur with the scaled variables, for example, 
using
\begin{eqnarray}   
  \sfrac{\partial \tilde{U}_{n+1}^F}{\partial \lambda} &=& 
  (0.5 / \sigma_n) \left[ \begin{array}{cc} 
       -h_n \phi_n(h_n) & (h_n \phi_{n,x}(h_n) - \phi_n(h_n)) / \tau_n \\ 
       -\phi_n(h_n) - h_n \phi_{n, x}(h_n) & -h_n \phi_n(h_n) 
        \end{array} \right]  \tilde{U}_n^F  \nonumber \\
    & & + (1 / \sigma_n) \left[ \begin{array}{cc} \phi_{n,x}(h_n) & \phi_n(h_n) \\ 
       -\tau_n \phi_n(h_n) & \phi_{n,x}(h_n) \end{array} \right] 
        \sfrac{\partial \tilde{U}_n^F}{\partial \lambda}  \label{SC_FWRD}
\end{eqnarray} 
instead of (\ref{VAR_FWRD_N}).  For the forward recurrences we take 
$n = 1, 2, \ldots, M-1$, while for the backward recurrences 
$n = N, N-1, \ldots, M$.   

It remains to recover the desired values 
\begin{eqnarray*}
  \left[ \begin{array}{c} u_{\lambda}(\ell) \\ u_{x,\lambda}(\ell) \end{array}  
         \right] &=& \sfrac{\partial U_{N+1}^F}{\partial \lambda} \\
  \left[ \begin{array}{c} v_{\lambda}(\ell) \\ v_{x,\lambda}(\ell) \end{array}  
         \right] &=& \sfrac{\partial V_{N+1}^F}{\partial \lambda}
\end{eqnarray*}
from the scaled variables.

From the first component of (\ref{LC1}) we have for any $x \in [0, \ell]$
\begin{eqnarray*}
 u_{\lambda} &=& c_{11} u_{\lambda}^B + c_{12} v_{\lambda}^B + 
       u^B \sfrac{\partial c_{11}}{\partial \lambda} 
     + v^B \sfrac{\partial c_{12}}{\partial \lambda} \\
     &=& \sfrac{\partial c_{11}}{\partial \lambda} \qquad \mbox{ at } x = \ell.
\end{eqnarray*}
Similarly
\begin{eqnarray*}
 u_{x \lambda} &=& c_{11} u_{x \lambda}^B + c_{12} v_{x \lambda}^B + 
       u_x^B \sfrac{\partial c_{11}}{\partial \lambda} 
     + v_x^B \sfrac{\partial c_{12}}{\partial \lambda}  \\
     &=& \sfrac{\partial c_{12}}{\partial \lambda} \qquad \mbox{ at } x = \ell,
\end{eqnarray*}
\begin{eqnarray*}
 v_{\lambda} &=& c_{21} u_{\lambda}^B + c_{22} v_{\lambda}^B + 
       u^B \sfrac{\partial c_{21}}{\partial \lambda} 
     + v^B \sfrac{\partial c_{22}}{\partial \lambda}  \\
     &=& \sfrac{\partial c_{21}}{\partial \lambda} \qquad \mbox{ at } x = \ell,
\end{eqnarray*}
and
\begin{eqnarray*}
 v_{x \lambda} &=& c_{21} u_{x \lambda}^B + c_{22} v_{x \lambda}^B + 
       u_x^B \sfrac{\partial c_{21}}{\partial \lambda} 
     + v_x^B \sfrac{\partial c_{22}}{\partial \lambda}   \\
     &=& \sfrac{\partial c_{22}}{\partial \lambda} \qquad \mbox{ at } x = \ell.
\end{eqnarray*}
Finally
\[   \sfrac{\partial c_{11}}{\partial \lambda} = -(v_x^B u^F - v^B u_x^F)
   \sfrac{\partial \Delta}{\partial \lambda} / \Delta^2 
   + [v_{x,\lambda}^B u^F + v_x^B u_{\lambda}^F 
   - v_{\lambda}^B u_x^F - v^B u_{x, \lambda}^F] / \Delta  \]
and 
\[   \sfrac{\partial \Delta}{\partial \lambda} = u_{\lambda}^B v_x^B 
  + u^B v_{x \lambda}^B - u_{x \lambda}^B v^B - u_x^B v_{\lambda}^B.  \]
Similarly,
\[ \sfrac{\partial c_{12}}{\partial \lambda} = -(u^B u_x^F - u_x^B u^F)
   \sfrac{\partial \Delta}{\partial \lambda} / \Delta^2 
   + [u_{\lambda}^B u_x^F + u^B u_{x \lambda}^F 
   - u_{x \lambda}^B u^F - u_x^B u_{\lambda}^F] / \Delta  \]
\[ \sfrac{\partial c_{21}}{\partial \lambda} = -(v_x^B v^F - v^B v_x^F)
   \sfrac{\partial \Delta}{\partial \lambda} / \Delta^2 
   + [v_{x \lambda}^B v^F + v_x^B v_{\lambda}^F 
   - v_{\lambda}^B v_x^F - v^B v_{x \lambda}^F] / \Delta  \]
and
\[ \sfrac{\partial c_{22}}{\partial \lambda} = -(u^B v_x^F - u_x^B v^F)
   \sfrac{\partial \Delta}{\partial \lambda} / \Delta^2 
   + [u_{\lambda}^B v_x^F + u^B v_{x \lambda}^F 
   - u_{x \lambda}^B v^F - u_x^B v_{\lambda}^F] / \Delta.   \]
These last five are all to be evaluated at $x = x_M$.  By inspecting the 
scale factors, it follows that, as in Section~2, we must multiply by 
$\zeta_M$ given by (\ref{ZETA}) when using the scaled variables.
  
After the double-shooting and matching with scaled variables, the results 
are to be substituted into (\ref{FIX_D}) or (\ref{FIX_N}).  While these
formulas seem complicated, for this paper they are only to be used in the 
neighborhood of a 0/0.  The computation of $\lambda^*$ itself may be done 
using a characterization in terms of eigenvalues (see \cite{TOMSEAST}), or 
by searching for zeros of the numerators in the expressions for $f$, 
(\ref{f4}), i.e., zeros of $2 - |u^F(\ell, \lambda) + 
{v^F}'(\ell, \lambda)| = 2 - |c_{11} + c_{22}|$.  

Since we have no test problems with closed form solutions to verify computer
output for the variational variables, we have compared our algorithm with 
finite difference approximations.  Table~7.1 contains data for several of 
our examples with a Dirichlet initial condition.  In all cases a central 
difference was used with a stepsize of $10^{-4}$, and an absolute error 
tolerance of $10^{-8}$ was used for $u_{x \lambda}$ and $v_{\lambda}$.  
The agreement is good in all cases.

\vskip 0.2in
\begin{center}
Table 7.1. Finite difference estimates compared to $u_{x \lambda}$ and $v_{\lambda}$.
\vskip 4pt

\begin{tabular}{cr@{}c@{}lr@{}c@{}lr@{}c@{}lr@{}c@{}lr@{}c@{}l}
\hline
Example&\multicolumn{3}{c}{$\lambda$}&\multicolumn{3}{c}{$\Delta_\lambda u_x$}&
\multicolumn{3}{c}{$u_{x \lambda}$}&\multicolumn{3}{c}{$\Delta_\lambda v$}&
\multicolumn{3}{c}{$v_{\lambda}$}\\
\hline
(\ref{MATHIEU})&$-$0&.&35 &$-$63&.&7915&$-$63&.&7916&$-$56&.&2402&$-$56&.&24019 \\
               &   1&.&00 &$-$1&.&6844 &$-$1&.&684311&$-$5&.&2775&   5&.&277455 \\
               &   2&.&00 &   1&.&3092 &   1&.&309169&$-$2&.&2702&$-$2&.&270148 \\
  (\ref{EX3})  &   2&.&00 &$-$1&.&7131 &$-$1&.&713098&   2&.&3125&   2&.&312439 \\
               &   3&.&00 &   0&.&9515 &  0&.&9514705&   0&.&0938&   0&.&009380 \\
               &   5&.&00 &$-$2&.&6099 &$-$2&.&609927&   0&.&6906&   0&.&690553 \\
  (\ref{EX5})  &$-$0&.&40 &$-$5&.&870  &$-$5&.&870013&$-$112&.&35&$-$112&.&3457 \\
               &   1&.&00 &   2&.&6079 &$-$2&.&607899&   5&.&5210&   5&.&521029 \\
               &   2&.&00 &   1&.&9781 &   1&.&978065&$-$1&.&8360&$-$1&.&836113
\end{tabular}
\end{center}
\vskip 0.2in

We conclude this section with numerical data illustrating the overhead 
required for the additional calculation of the variational variables 
$u_{\lambda}^F$, $u_{x \lambda}^F$, $v_{\lambda}^F$ and $v_{x \lambda}^F$.  
Table~7.2 has timing data for all five examples where 601 $f(\lambda)$ 
evaluations were made over a uniform grid of $\lambda$ in the intervals 
shown.  Examples (\ref{MATHIEU}), (\ref{EX3}), and (\ref{EX5}) had Dirichlet
initial conditions; the other two had Neumann.  The column labelled
`basic' gives the time required for just the $u^F$, $u_x^F$, $v^F$, 
and $v_x^F$ calculations; the final column gives the time required to 
compute all eight variables.           

\pagebreak
\vskip 0.2in
\begin{center}
Table 7.2. Timings for calculation of basic and variational solutions.
\vskip 4pt

\begin{tabular}{rcccc}
\hline
    Example    &   Interval  & Tolerance &   Basic  &    All   \\
\hline
(\ref{MATHIEU})&  [1.2, 7.2] & $10^{-4}$ &   0.027  &   0.031  \\
               &             & $10^{-6}$ &   0.048  &   0.086  \\
               &             & $10^{-8}$ &   0.103  &   0.120  \\
               &             & $10^{-10}$&   0.126  &   0.149  \\
  (\ref{EX2})  &  [3.0, 9.0] & $10^{-4}$ &   0.036  &   0.041  \\
               &             & $10^{-6}$ &   0.078  &   0.121  \\
               &             & $10^{-8}$ &   0.146  &   0.164  \\
               &             & $10^{-10}$&   0.180  &   0.203  \\
  (\ref{EX3})  & [2.0, 8.0]  & $10^{-4}$ &   0.040  &   0.046  \\
               &             & $10^{-6}$ &   0.090  &   0.117  \\
               &             & $10^{-8}$ &   0.162  &   0.178  \\
               &             & $10^{-10}$&   0.198  &   0.222  \\
  (\ref{EX4})  & [1.2, 7.2]  & $10^{-4}$ &   0.048  &   0.054  \\
               &             & $10^{-6}$ &   0.143  &   0.162  \\
               &             & $10^{-8}$ &   0.196  &   0.214  \\
               &             & $10^{-10}$&   0.368  &   0.524  \\
  (\ref{EX5})  &  [1.5, 7.5] & $10^{-4}$ &   0.050  &   0.056  \\
               &             & $10^{-6}$ &   0.138  &   0.169  \\
               &             & $10^{-8}$ &   0.204  &   0.224  \\
               &             & $10^{-10}$&   0.336  &   0.510  \\
\hline
  Totals       &             & $10^{-4}$ &   0.201  &   0.228  \\
               &             & $10^{-6}$ &   0.497  &   0.655  \\
               &             & $10^{-8}$ &   0.811  &   0.900  \\
               &             & $10^{-10}$&   1.208  &   1.608  \\
\end{tabular}
\end{center}
\vskip 0.2in

Despite the doubling in the number of dependent variables, the overhead 
is increased by only 10\% to 30\% in the totals, largely because the basic
and variational variables use the same transcendental function values for 
a fixed $\lambda$.

\end{document}